\documentclass{amsart}

\usepackage{style}

\def\Enr{\mathrm{Enr}}

\title{Enriques involutions on pencils of K3 surfaces}

\author{Dino Festi}
\address[Dino Festi]
{Dipartimento di matematica Federigo Enriques, Università degli Studi di Milano, via Saldini 50, 20133 Milan, Italy}
\email{dino.festi@unimi.it}

\author{Davide Cesare Veniani}
\address[Davide Cesare Veniani]
{Institut für Topologie und Geometrie, 
Universität Stuttgart,
Pfaffenwaldring~57,
70569 Stuttgart, Germany}
\email{davide.veniani@mathematik.uni-stuttgart.de}

\date{\today}

\makeatletter
\@namedef{subjclassname@2020}{%
 \textup{2020} Mathematics Subject Classification}
\makeatother
\subjclass[2020]{%
14J28 % K3 surfaces and Enriques surfaces
(%
14J27 %Elliptic surfaces, elliptic or Calabi-Yau fibrations
14J50% Automorphisms of surfaces and higher-dimensional varieties
)}

\keywords{K3 surface, Enriques surface, transcendental lattice, elliptic fibration}

\begin{document}

\maketitle

\begin{abstract}
The three pencils of K3 surfaces of minimal discriminant whose general element covers at least one Enriques surface are Kond\={o}'s pencils I and II, and the Apéry--Fermi pencil. We enumerate and investigate all Enriques surfaces covered by their general elements. 
\end{abstract}

%\tableofcontents

\section{Introduction}

Any complex Enriques surface is doubly covered by a K3 surface. 
On the other hand, a K3 surface \(X\) can cover infinitely many Enriques surfaces. The set \(\Enr(X)\) of isomorphism classes of Enriques surfaces doubly covered by \(X\), though, is always finite by a result of Ohashi \cite{Ohashi:number.Enriques}.
We call its cardinality \(|{\Enr(X)}|\) the \emph{Enriques number} of the K3 surface \(X\).
The Enriques number \(|{\Enr(X)}|\) only depends on the transcendental lattice of~\(X\). Shimada and the second author~\cite{Shimada.Veniani:Enriques.involutions.singular.K3} described a procedure to determine \(|{\Enr(X)}|\) and applied it to K3 surfaces of maximal Picard rank~\(20\).

A K3 surface \(X\) of Picard rank~\(19\) can be seen as the generic element of a pencil of K3 surfaces. Its transcendental lattice \(T_X\) is an even lattice of signature \((2,1)\).
By a result of Brandhorst, Sonel and the second author~\cite{Brandhorst.Sonel.Veniani:idoneal.genera}, the surface \(X\) covers an Enriques surface only if \(4\) divides \(\det(T_X)\), but this condition is not sufficient. In this paper we analyze in detail what happens when \(|{\det(T_X)}|\) is small, more precisely 
\begin{equation} \label{eq:det(T_X)<16}
    |{\det(T_X)}| < 16.
\end{equation}

Henceforth, let \(X\) be a K3 surface of Picard rank \(19\) with transcendental lattice \(T_X\). In the case \(T_X \cong \bU \oplus [2n]\), \(n \geq 1\), it was already noted by Hulek and Schütt~\cite{Hulek.Schuett:Enriques.surf.jacobian.ell} that \(\Enr(X) \neq \emptyset\) if and only if \(n\) is even.
Indeed, we prove in \autoref{lem:|det(T_X)|<16} under assumption~\eqref{eq:det(T_X)<16} that \(\Enr(X) \neq \emptyset\) if and only if 
\[
    T_X \cong \bU \oplus [4],\,\bU \oplus [8] \text{ or } \bU \oplus [12].
\]

The main reason for bound \eqref{eq:det(T_X)<16} is to keep computations feasible.
In particular the enumeration of jacobian elliptic fibrations on K3 surfaces with \(T_X \cong \bU \oplus [16]\) already becomes quite hard. Moreover, the pencil of K3 surfaces with \(T_X \cong \bU(2) \oplus [4]\) is not of the form \(\bU \oplus [2n]\), but it still holds \(\Enr(X) \neq \emptyset\), as its generic element is a Kummer surface~\cite{Keum:Kummer.Enriques}.

Quite interestingly, the first two pencils already feature prominently in Kond\={o}'s classification of Enriques surfaces with finite automorphism group \cite{Kondo:Enriques.finite.aut}, which we briefly recall.
There are seven families of such Enriques surfaces, numbered \(\I\) to \(\VII\). 
Families \(\I\) and \(\II\) are \(1\)-dimensional, while families \(\III\) to \(\VII\) are \(0\)-dimensional. 
The K3 surfaces covering the generic Enriques surface of type~\(\I\) and \(\II\) have transcendental lattice \(T_X \cong \bU \oplus [4]\) and \(T_X \cong \bU \oplus [8]\), respectively. 

The third pencil with \(T_X \cong \bU \oplus [12]\) has also been extensively studied, because of its arithmetical properties and its appearance in several seemingly unrelated physical contexts (see \cite{Festi.vanStraten:Bhabha.scattering, Peters.Stienstra:pencil.Apery}). 
Following Bertin and Lecacheux~\cite{Bertin.Lecacheux:apery-fermi.2-isogenies}, who classified the elliptic fibrations on its generic element (\autoref{tab:genus_Apery_Fermi}), we call it the \emph{Apéry--Fermi pencil}.

The aim of this paper is to enumerate and investigate the Enriques surfaces covered by these three pencils. 
More precisely, for each \(m \in \set{1,2,3}\) we consider a K3 surface \(X\) with \(T_X \cong \bU \oplus [4m]\) and do the following:
\begin{itemize}
    \item we compute the Enriques number \(|{\Enr(X)}|\);
    \item we classify all jacobian elliptic fibrations on \(X\) using the extension of the Kneser--Nishiyama method explained in \cite{Festi.Veniani:counting.ell.fibr.K3};
    \item we relate the special elliptic pencils on the Enriques quotients to the elliptic fibrations on~\(X\).
\end{itemize}

We summarize here our findings.

Fix \(m \in \IZ\), \(m \geq 1\), and let \(\omega\) be the number of prime divisors of \(2m\) and \(X\) a K3 surface with \(T_X \cong \bU \oplus [4m]\), \(m \geq 1\).
Among the Enriques quotients of \(X\) there are \(2^{\omega-1}\) which we call \emph{of Barth--Peters type} (\autoref{lem:BP.quotients}).
Such quotients admit a cohomologically trivial involution (see~\cite{Mukai:numerically.trivial.inv.Kummer.type,Mukai.Namikawa:aut.Enriques.which.act.trivially}) and their presence is explained by the fact that our pencils are subfamilies of the \(2\)-dimensional Barth--Peters family, a fact already noted by Hulek and Schütt \cite{Hulek.Schuett:Enriques.surf.jacobian.ell,Hulek.Schuett:arithmetic.singular.Enriques}.

It turns out that if \(m = 1\), then \(X\) covers only one Enriques surface \(Y\) (\autoref{thm:Enr(X)_Kondo_I}). Therefore, the Enriques surface \(Y\) is of Barth--Peters type and, moreover, coincides with Kond\={o}'s quotient, so it has finite automorphism group.
The list of the \(9\) elliptic fibrations on \(X\) appears in other papers by Scattone~\cite{Scattone:on.the.compactification}, Dolgachev~\cite{Dolgachev:Mirror.symm.latt.pol.K3} and Elkies and Schütt~\cite{Elkies.Schuett:K3.families.high.Picard}, and we confirm it here (\autoref{tab:genus_Kondo_I}).

If \(m = 2\), then \(X\) covers two Enriques surfaces \(Y',Y''\), of which only one, say \(Y'\), is of Barth--Peters type. We show that the other surface \(Y''\) is Kond\={o}'s quotient with finite automorphism group. 
We include the classification of elliptic fibrations on \(X\) up to automorphisms (\autoref{tab:genus_Kondo_II}).
One subtlety arises in this case: two of the \(17\) elliptic fibrations on \(X\) (No. 12 and 13 in \autoref{tab:genus_Kondo_II}) have the same Mordell-Weil group and two singular fibers of type~\(\I_4^*\). 
Nonetheless, the two fibrations are not equivalent under the action of \(\Aut(X)\), as they have different frames.
We determine which one is the pullback of a special elliptic pencil on \(Y'\) and which one is the pullback of a special elliptic pencil on \(Y''\) (\autoref{rmk:W11,W12}).

Finally, if \(m = 3\), then \(X\) covers three Enriques surfaces \(Y',Y'',Y'''\), of which two, say \(Y'\) and \(Y''\) are of Barth--Peters type (\autoref{thm:Apery_Fermi}).
Applying a construction by Hulek and Schütt \cite[§3]{Hulek.Schuett:Enriques.surf.jacobian.ell} and using a particular configuration of curves on \(X\) found by Peters and Stienstra, we determine a simple description of an explicit Enriques involution for \(Y'''\). 
In this way we find a configuration of smooth rational curves on \(Y'''\) whose dual graph is the union of a tetrahedron and a complete graph of degree~\(6\) (\autoref{rmk:Apery-Fermi.ell.fibr}).

\subsection*{Acknowledgments}
We warmly thank Simon Brandhorst, Klaus Hulek, Matthias Schütt and Ichiro Shimada for their valuable comments. We are also grateful to the anonymous referee for carefully reading the manuscript and for their useful remarks.

\section{Preliminary results}
In this section, after explaining our conventions on lattices in \autoref{sec:lattices}, we collect results regarding K3 surfaces with transcendental lattice \(T_X \cong \bU \oplus [2m]\), \(m \in \IZ\), especially regarding their jacobian elliptic fibrations in \autoref{sec:fibrations}.
In~\autoref{sec:Enriques.numbers} we recall the enumeration formula for Enriques quotients contained in \cite{Shimada.Veniani:Enriques.involutions.singular.K3} and we prove the lemma that motivates the whole paper.
Finally, in \autoref{sec:Barth-Peters.type} we introduce the notion of Enriques quotient of Barth--Peters type.

\subsection{Lattices} \label{sec:lattices}
In this paper, a \emph{lattice of rank \(r\)} is a finitely generated free \(\IZ\)-module \(L \cong \IZ^r\) endowed with an integral symmetric bilinear form \(L \times L \rightarrow \IZ\) denoted \((v,w) \mapsto v \cdot w\). 
The \emph{signature} of \(L\) is the signature of the induced real symmetric form on \(L \otimes \IR\).
We say that \(L\) is \emph{even} if \(v^2 \coloneqq v \cdot v \in 2 \IZ\) for every \(v \in L\).
The \emph{dual} \(L^\vee \coloneqq \hom(L,\IZ)\) of \(L\) can be identified with \(\Set{w \in L \otimes \IQ}{w \cdot v \in \IZ \text{ for all } v \in L}\). The \emph{discriminant group} of \(L\) is defined as
\[
    L^\sharp \coloneqq L^\vee/L,
\]
which is a finite abelian group. 
We denote by \(\ell(L^\sharp)\) its \emph{length}, i.e. the minimal number of generators. For a prime number \(p\) we denote by \(\ell_p(L^\sharp)\) its \emph{\(p\)-length}, i.e. the minimal number of generators of its \(p\)-part.

If \(L\) is an even lattice, then \(L^\sharp\) acquires naturally the structure of a finite quadratic form \(L^\sharp \rightarrow \IQ/2\IZ\). 
There is a natural homomorphism \(\OO(L) \rightarrow \OO(L^\sharp)\) denoted \(\gamma \mapsto \gamma^\sharp\). 

We write \(\bU\) for the indefinite unimodular even lattice of rank~\(2\), and \(\bA_n,\bD_n,\bE_n\) for the negative definite ADE lattices. The notation~\([m]\), with \(m \in \IZ\), denotes the lattice of rank~\(1\) generated by a vector of square~\(m\).
We adopt Miranda--Morrison's notation~\cite{Miranda.Morrison} for the elementary finite quadratic forms \(\bu_k,\bv_k,\bw^\varepsilon_{p,k}\). 
We recall that \(\bu_k\) (resp. \(\bv_k\)) is generated by two elements of order \(2^k\), both of square \(0 \in \IQ/2\IZ\) (resp. \(1 \in \IQ/2\IZ\)), such that their product is equal to \(1/2^k \in \IQ/\IZ\). The forms \(\bw_{2,k}^\varepsilon\), with \(\varepsilon \in \set{1,3,5,7}\), are generated by one element of order \(2^k\) and square \(\varepsilon/2^k \in \IQ/2\IZ\). 
For an odd prime \(p\) the forms \(\bw_{p,k}^\varepsilon\), with \(\varepsilon \in \set{\pm 1}\), are generated by one element of order \(p^k\) and square \(a/p^k \in \IQ/2\IZ\), where \(a\) is a square modulo \(p\) if and only if \(\varepsilon = 1\).

The \emph{genus} of a lattice \(L\) is defined as the set of isomorphism classes of lattices \(M\) with \(\sign(L) = \sign(M)\) and \(L^\sharp \cong M^\sharp\). 
A genus is always a finite set (see \cite[Satz~21.3]{Kneser:quadratische.Formen}).

An embedding of lattices \(\iota\colon  M \hookrightarrow L\) is called \emph{primitive} if \(L/\iota(M)\) is a free group. We denote by \(\iota(M)^\perp \subset L\) the \emph{orthogonal complement} of \( M \) inside \( L \). 
We quickly summarize Nikulin's theory of primitive embeddings \cite{Nikulin:int.sym.bilinear.forms}. 

By \cite[Prop.~1.5.1]{Nikulin:int.sym.bilinear.forms} a primitive embedding of even lattices \(M \hookrightarrow L\) is given by a subgroup \(H \subset M^\sharp\) and an isometry \[
\gamma \colon H \rightarrow H' \coloneqq \gamma(H) \subset (\iota(M)^\perp(-1))^\sharp.
\]
If \(\Gamma\) denotes the graph of \(\gamma\) in \(M^\sharp \oplus (\iota(M)^\perp(-1))^\sharp\), the following identification between finite quadratic forms holds (the finite quadratic form on the right side being induced by the one on \(M^\sharp \oplus (\iota(M)^\perp(-1))^\sharp\)):
\begin{equation} \label{eq:L^sharp=Gamma^perp/Gamma}
  L^\sharp \cong \Gamma^\perp/\Gamma.  
\end{equation}
In this paper we call \(H\), \(\gamma\) resp. \(\Gamma\) the \emph{gluing subgroup}, \emph{gluing isometry} resp. \emph{gluing graph} of \(M \hookrightarrow L\).

Equivalently by \cite[Prop.~1.15.1]{Nikulin:int.sym.bilinear.forms}, assuming that \(L\) is unique in its genus, a primitive embedding \(M \hookrightarrow L\) is given by a subgroup \(K \subset L^\sharp\) and an isometry 
\[
    \xi\colon K \rightarrow K' \coloneqq \xi(K) \subset M(-1)^\sharp.
\]
If \(\Xi\) denotes the graph of \(\xi\) in \(L^\sharp \oplus M(-1)^\sharp\), the following identification between finite quadratic forms holds (the finite quadratic form on the right side being induced by the one on \(L^\sharp \oplus M(-1)^\sharp\)):
\begin{equation} \label{eq:(M^perp)^sharp=Xi^perp/Xi}
    (\iota(M)^\perp)^\sharp \cong \Xi^\perp/\Xi.
\end{equation}
In this paper we call \(K\), \(\xi\) resp. \(\Xi\) the \emph{embedding subgroup}, \emph{embedding isometry} resp. \emph{embedding graph} of \(M \hookrightarrow L\).

\subsection{Elliptic fibrations} \label{sec:fibrations}

Given a K3 surface \(X\), we denote \(T_X\) its transcendental lattice, \(S_X\) its Néron--Severi lattice, and \(\cJ_X\) the set of jacobian elliptic fibrations on \(X\).
The \emph{frame genus} of \(X\) is defined as the genus \(\cW_X\) of negative definite lattices~\(W\) with \(\rank(W) = \rank(S_X) - 2\) and \(W^\sharp \cong S_X^\sharp\). 
The lattices in \(\cW_X\) are called \emph{frames}.
The classes of a fiber and a section of a jacobian elliptic fibration induces a primitive embedding \(\iota \colon \bU \hookrightarrow S_X\).
As explained in \cite{Festi.Veniani:counting.ell.fibr.K3}, there is a well-defined function
\[
    \fr_X \colon \cJ_X/{\Aut(X)} \rightarrow \cW_X
\]
which sends each jacobian fibration to the isomorphism class of \(\iota(\bU)^\perp \subset S_X\).

\begin{lemma} \label{lem:BP.fibrations}
If \(X\) is a K3 surface with transcendental lattice \(T_X \cong \bU \oplus [2n]\), \(n \geq 1\),
then on \(X\) there are exactly \(2^{\omega-1}\) jacobian elliptic fibrations with frame \(W \coloneqq \bE_8^2 \oplus [-2n]\) up to automorphisms, where \(\omega\) is the number of prime divisors of~\(2n\).
\end{lemma}
\proof
Essentially by \cite[Thm.~2.8]{Festi.Veniani:counting.ell.fibr.K3} we want to prove that
\[
    |{\OO_\hdg^\sharp(T_X)}\backslash{\OO(T_X^\sharp)}/{\OO(W)}| = 2^{\omega-1},
\]

If \(n = 1\), then \(\OO(T_X^\sharp) = \set{\id}\) and we conclude immediately.

Suppose that \(n \geq 2\). Since \(\ell(T_X^\sharp) = 1\), the discriminant form \(T_X^\sharp\) is the direct sum of forms \(\bw_{p,k}^\varepsilon\).
It holds \(\OO(q \oplus q') \cong \OO(q) \times \OO(q')\) if \(q\) and \(q'\) are finite quadratic forms with \(|q|\) and \(|q'|\) coprime, and \(|{\OO(\bw^{\varepsilon}_{p,k})}| = 2\) if \(p\) is odd or \(p = 2\) and \(k \geq 2\).
Hence, \(\OO(T_X^\sharp)\) is a \(2\)-elementary group of length~\(\omega\). In particular,
\[
    |{\OO(T_X^\sharp)}| = 2^{\omega},
\]
As \(\rank(T_X)\) is odd, it holds \(\OO_\hodge^\sharp(T_X) = \set{\pm \id}\) (see for instance \cite[Cor.~3.3.5]{Huybrechts:lectures.K3}). 
Note, moreover, that \(\id \neq -\id\) in \(T_X^\sharp\).
The orthogonal group of \(W\) is the direct sum of \(\OO(\bE_8^2)\), which has trivial action on the discriminant group because \(\bE_8\) is unimodular, and \(\OO([-2n]) = \set{\pm \id}\). 
Therefore, it also holds \(\OO^\sharp(W) = \set{\pm \id}\), so we have
\[
    |{\OO_\hodge^\sharp(T_X)}\backslash{\OO(T_X^\sharp)}/{\OO^\sharp(W)}| = |{\OO(T_X^\sharp)}/{\set{\pm \id}}| = |{\OO(T_X^\sharp)}| / |{\set{\pm \id}}| = 2^{\omega-1}. \qedhere
\]
\endproof 

The \emph{Mordell--Weil group}, i.e. the group of sections of a jacobian elliptic fibration, is naturally endowed with a rational symmetric bilinear form denoted by \((P,Q) \mapsto \langle P, Q \rangle \in \IQ\), called the \emph{Mordell--Weil lattice}. The \emph{height} of a section is defined as \(\height(P) \coloneqq \langle P, P \rangle\). For a clear exposition of this topic we refer to Shioda's original paper~\cite{Shioda:on.MW.lattices}.

\begin{remark} \label{rmk:Hulek.Schutt.II^*}
Consider one of the elliptic fibrations \(\pi \colon X \rightarrow \IP^1\) as in \autoref{lem:BP.fibrations} and for simplicity assume that \(n \geq 2\). 
Since \(W_\rootlattice \cong \bE_8^2\), the fibration \(\pi\) has two singular fibers of Kodaira type~\(\II^*\).
As already remarked by Hulek and Schütt~\cite[§4.2.2]{Hulek.Schuett:Enriques.surf.jacobian.ell}, starting from the fibration \(\pi\) we can construct an involution on \(X\) which turns out to be an Enriques involution if \(n\) is even. 
We repeat here their construction directly on the lattice \(S_X \cong \bU \oplus \bE_8^2 \oplus [-2n]\).
In the following computations we let \(\OO(S_X)\) act on \(S_X\) \emph{from the right}, 
so the composition of two isometries in \(\OO(S_X)\) corresponds to the product of their associated matrices in reversed order.

Let \(s_1,\ldots,s_{19}\) be a system of generators of \(S_X\) such that the corresponding Gram matrix is the standard one.
Then, \(S_X^\sharp\) is generated by \(s_{19}/(2n)\).
In these coordinates, consider the vectors
\begin{align*} 
    F & \coloneqq ( 1, 0, 0,0,0,0,0,0,0,0,0,0,0,0,0,0,0,0,0), \\
    O & \coloneqq (-1, 1,0,0,0,0,0,0,0,0,0,0,0,0,0,0,0,0,0), \\
    P & \coloneqq (n-1,1,0,0,0,0,0,0,0,0,0,0,0,0,0,0,0,0,1).
\end{align*}
Note that \(F^2 = 0, O^2 = P^2 = -2, F \cdot O = F \cdot P = 1, P \cdot O = n-2\).

We can suppose that \(s_3,\ldots,s_{18}\), generating the two copies of \(\bE_8\), correspond to the components of the singular fibers which do not intersect \(O\), \(F\) to the class of a fiber, \(O\) to a section which we take as origin and \(P\) to a section of height (cf. \cite[eq. (8.19)]{Shioda:on.MW.lattices})
\[
 \height(P) = \langle P,P \rangle = 2 \chi(\cO_X) + 2 P \cdot O = 2 \cdot 2 + 2 \cdot (n-2) = 2n.
\]

Let \(Q \coloneqq \boxminus P\) be the inverse section of \(P\). Then \(\height(Q) = \height(P) = 2n\), hence \(Q\cdot O = n-2\). 
It follows from \(\langle P,Q \rangle = - \langle P,P \rangle\) that (cf. \cite[eq. (8.18)]{Shioda:on.MW.lattices})
\[
    P \cdot Q = \chi(\cO_X) + P \cdot O + Q \cdot O - \langle P, Q \rangle = 2 + (n-2) + (n-2) + 2n = 4n- 2.
\]
Recalling moreover that \(F \cdot Q = 1\) and \(Q\cdot Q = -2\), we see that in our basis we can write
\[
    Q = (n-1,1,0,0,0,0,0,0,0,0,0,0,0,0,0,0,-1).
\]
Let \(t_{P} \in \OO(S_X)\) be the pullback of the automorphism induced by the translation by \(P\). We have \(t_P(Q) = O, t_P(O) = P, t_P(F) = F\) and \(t_P\) acts trivially on the other components of the fibers of type~\(\II^*\). 
Therefore, \(\tau_P\) is given by the following matrix:
\[
t_P = \begin{pmatrix}
    1 &   & &  \\
    n & 1 & & 1 \\
    & & I_{16} & \\
    2n &   & & 1 \\
\end{pmatrix}
\]

Let now \(\imath \in \OO(S_X)\) be the isometry given by
\[
\imath = \begin{pmatrix}
    I_2 & \\
    & & I_{8} & \\
    & I_8 & \\
    & & & -1 \\
\end{pmatrix}.
\]
Note that \(\imath\) swaps the two fibers of type \(\II^*\) and that \(\imath(F) = F\), \(\imath(P) = Q\) and \(\imath(Q) = P\). 
Constructing an ample divisor as in \cite[Prop.~2.7]{Festi.Veniani:counting.ell.fibr.K3}, it is easy to see that \(\imath\) preserves the ample cone. Moreover, \(\imath^\sharp = -\id \in S_X^\sharp\). Therefore, by the Torelli theorem \(\imath\) is the pullback of a non-symplectic involution (whose quotient is a rational surface).

Consider \(\varepsilon \coloneqq \imath \circ t_P\), whose matrix is then given by 
\[
\varepsilon = \begin{pmatrix}
    1 \\
    n & 1 & & & -1\\
    & & & I_{8} & \\
    & & I_8 & \\
    2n & & & & -1 \\
\end{pmatrix}.
\]
A computation shows that for even \(n\), the invariant lattice of \(\varepsilon\) is isometric to \(\bE_8(2) \oplus \bU(2)\), hence it corresponds to an Enriques involution by Nikulin's classification~\cite[Thm.~4.2.2]{Nikulin:quotient.groups.hyperbolic.forms.2-reflections}.
% (see the file \verb+calc/iota_tauP.sage+).
The coinvariant lattice is isometric to \(\bE_8(2) \oplus [-2n]\).
\end{remark}

\subsection{Enriques numbers} \label{sec:Enriques.numbers}
Let \(X\) be a K3 surface with Néron--Severi lattice \(S_X\) and transcendental lattice \(T_X\).
We recall briefly the formula for \(|{\Enr(X)}|\) proved in \cite{Shimada.Veniani:Enriques.involutions.singular.K3}. 
We define
\[
    \bM \coloneqq \bU(2) \oplus \bE_8(2).
\]

By Nikulin's classification~\cite{Nikulin:quotient.groups.hyperbolic.forms.2-reflections}, if \(\varepsilon \in \OO(S_X)\) is the pullback of an Enriques involution, then the invariant sublattice 
\(S_X^\varepsilon \coloneqq \Set{x \in S_X}{\varepsilon(x) = x}\) 
is isomorphic to \(\bM\). We denote by \((S_X)_\varepsilon \coloneqq (S_X^\varepsilon)^\perp\) the coinvariant lattice, whose isometry class depends on the involution \(\varepsilon\).

Given a primitive embedding \(\iota \colon \bM \hookrightarrow S_X\), we put
\[
    \OO(S_X,\iota) \coloneqq \{\, \varphi \in \OO(S_X) \,\vert\, \varphi(\iota(\bM)) = \iota(\bM)\,\}.
\]
The Hodge structure on \(\HH^2(X,\IZ)\) induces a Hodge structure on \(T_X\). We write \(\OO_\hodge(T_X)\) for the group of Hodge isometries of \(T_X\).
We fix an anti-isometry \(T_X^\sharp \cong S_X^\sharp\) (cf. \cite[Prop.~1.6.1]{Nikulin:int.sym.bilinear.forms}), so that we can identify \(\OO(T_X^\sharp) \cong \OO(S_X^\sharp)\). 
We denote the images of \(\OO_\hodge(T_X)\) and \(\OO(S_X,\iota)\) under the natural morphisms \(\OO(T_X) \rightarrow \OO(T_X^\sharp)\) and \(\OO(S_X) \rightarrow \OO(S_X^\sharp)\) by \(\OO_\hodge^\sharp(T_X)\) and \(\OO^\sharp(S_X,\iota)\), respectively.

\begin{theorem}[{\cite[Thm.~3.1.9]{Shimada.Veniani:Enriques.involutions.singular.K3}}] \label{thm:Shimada.Veniani}
For any K3 surface \(X\) it holds
\[
    |{\Enr(X)}| = \sum 
    |{\OO_\hodge^\sharp(T_X)}\backslash{\OO(T_X^\sharp)}/{\OO^\sharp(S_X,\iota)}|,
\]
where the sum runs over all primitive embeddings \(\iota\colon \bM \hookrightarrow S_X\) up to the action of \(\OO(S_X)\) such that there exists no \(v \in \iota(\bM)^\perp\) with \(v^2 = -2\).
\end{theorem}

The main topic of the present paper are K3 surfaces of Picard rank~\(19\) covering an Enriques surface. 
For computational reasons, we restrict ourselves to K3 surfaces whose transcendental lattice has discriminant \(|{\det(T_X)}| < 16\). By the next lemma, we have three cases to consider.

\begin{lemma} \label{lem:|det(T_X)|<16}
Let \(X\) be a K3 surface of Picard rank~\(19\) and transcendental lattice \(T_X\) and suppose that \(|{\det(T_X)}| < 16\). Then, \(\Enr(X) \neq \emptyset\) if and only if \(T_X \cong \bU \oplus [4m]\) with \(m \in \{1,2,3\}\).
\end{lemma}
\proof
If \(T_X \cong \bU \oplus [4m]\), \(m \geq 1\), then \(X\) covers an Enriques surface by~\cite[Proposition 4.2]{Hulek.Schuett:Enriques.surf.jacobian.ell}.

Conversely, suppose that \(\Enr(X) \neq \emptyset\). By~\cite[Thm.~1.1]{Brandhorst.Sonel.Veniani:idoneal.genera} the lattice \(T_X\) has a Gram matrix of the form
\[
    \begin{pmatrix}
        2a_{11} & a_{12} & a_{13} \\
        a_{12} & 4a_{22} & 2a_{23} \\
        a_{13} & 2a_{23} & 4a_{33},
    \end{pmatrix}, \quad a_{ij} \in \IZ.
\]
Therefore, \(\det(T_X)\) is divisible by \(4\).
Now, the discriminant group \(T_X^\sharp\) is a finite quadratic form on an abelian group of order \(|{\det(T_X)}|\) and of signature \(2-1 = 1\), because \(T_X\) has signature \((2,1)\). 
We classify such finite quadratic forms \(q\) using Miranda and Morrison's normal form \cite{Miranda.Morrison}. 
For each \(q\) in the list we find a lattice \(T\) such that \(T^\sharp \cong q\), obtaining the following table.
\begin{table}
    \centering
    \begin{tabular}{c|llllll}
    \toprule
    \({\det(T)}\) & \(-4\) & \(-8\) & \(-8\) & \(-8\) & \(-12\) & \(-12\) \\[2mm]
    \(T^\sharp\) & \(\bw^1_{2,2}\) & \(\bu_1 \oplus \bw^1_{2,1}\) & \(\bw^1_{2,3}\) & \(\bw^5_{2,3}\)& \(\bw^{3}_{2,2} \oplus \bw^{1}_{3,1}\) & \(\bw^{7}_{2,2} \oplus \bw^{-1}_{3,1}\) \\[2mm]
    \(T\) & \(\bU \oplus [4]\) & \(\bU(2) \oplus [2]\) & \(\bU \oplus [8]\) & {\small \(\begin{pmatrix} 2 & -1 & 0 \\ -1 & 2 & -1 \\ 0 & -1 & -2 \end{pmatrix}\)} & \(\bU \oplus [12]\)  & {\small \(\begin{pmatrix} 2 & -1 & 0 \\ -1 & -2 & -1 \\ 0 & -1 & 2 \end{pmatrix}\)} \\
    \bottomrule
    \end{tabular}
\end{table}

In all cases except the second one, \(T\) is unique by \cite[Thm.~1.14.2]{Nikulin:int.sym.bilinear.forms}. In the case \(T = \bU(2) \oplus [2]\), \(T\) is unique because \(T = T'(2)\), with \(T'\) a unimodular indefinite lattice.

If \(T_X \cong T\) is such that \(T^\sharp \cong \bw^5_{2,3}\) or \(T^\sharp \cong \bw^7_{2,2} \oplus \bw^{-1}_{3,1}\), then \(\Enr(X) = \emptyset\) because of \cite[Prop.~3.9]{Brandhorst.Sonel.Veniani:idoneal.genera} (in the notation of \cite{Brandhorst.Sonel.Veniani:idoneal.genera}, the two forms do not satisfy condition C(1)). The case \(T_X \cong \bU(2) \oplus [2]\) is excluded because of \cite[Thm.~1.1]{Brandhorst.Sonel.Veniani:idoneal.genera} (the lattice is an `exceptional lattice').

Therefore, the only cases left are \(T_X \cong T \in \set{\bU \oplus [4],\bU \oplus [8],\bU \oplus [12]}\).
\endproof

\subsection{Enriques quotients of Barth--Peters type} \label{sec:Barth-Peters.type}
Let now \(X\) be a K3 surface with \(T_X \cong \bU \oplus [4m]\), \(m \geq 1\).
A primitive embedding \(\iota \colon \bM \hookrightarrow S_X\) depends in general on several data (cf. \autoref{sec:lattices} and \cite[Prop.~1.15.1]{Nikulin:int.sym.bilinear.forms}), 
but in this case one only has to consider the orthogonal complement of the image, thanks to the next lemma.

\begin{lemma} \label{lem:emb.M-->S_X.for.T_X=U+4m}
Let \(X\) be a K3 surface with transcendental lattice \(T_X \cong \bU \oplus [4m]\), \(m \geq 1\).
If \(\iota \colon \bM \hookrightarrow S_X\) is a primitive embedding, then \(\iota(\bM)^\perp \cong N(2)\), where \(N\) is a lattice in the genus of \(\bE_8 \oplus [-2m]\). Conversely, for each such lattice \(N\) there exists exactly one primitive embedding \(\iota \colon \bM \hookrightarrow S_X\) with \(\iota(\bM)^\perp \cong N(2)\) up to the action of \(\OO(S_X)\). 
\end{lemma}
\proof
The Neron--Séveri lattice \(S_X\) is isomorphic to \(\bU \oplus \bE_8^2 \oplus [-4m]\).
Consider a primitive embedding \(\iota\colon \bM \hookrightarrow S_X\) with embedding subgroup \(K \subset S_X^\sharp\) and embedding graph \(\Xi\) (see~\autoref{sec:lattices}). 

Since \(\bM^\sharp \cong 5\bu_1\) is \(2\)-elementary and \(S_X^\sharp\) has length \(1\), it holds either \(|K| = 1\) or \(|K| = 2\). 
The first case, though, is impossible, as otherwise \((\iota(\bM)^\perp)^\sharp\) would have length \(11 > \rank(\iota(\bM)^\perp)\). 
Therefore, it must be \(|K| = 2\), so there is only one choice for the subgroup \(K \subset S_X^\sharp\), which is generated by an element of order \(2\) and square \(0 \in \IQ/2\IZ\).

Moreover, when taking \(\Xi^\perp/\Xi\) in the identification~\eqref{eq:(M^perp)^sharp=Xi^perp/Xi} one copy of \(\bu_1\) gets killed. Hence, it holds
\[
    (\iota(\bM)^\perp)^\sharp \cong 4\bu_1 \oplus [-4m]^\sharp,
\]
and in particular \(\ell_2((\iota(\bM)^\perp)^\sharp) = 9 = \rank(\iota(\bM)^\perp)\). 
Therefore (see for instance \cite[Lemma~3.10]{Brandhorst.Sonel.Veniani:idoneal.genera}), it holds \(\iota(\bM)^\perp \cong N(2)\), with \(N\) an even lattice. The genus of \(N(2)\) determines the genus of \(N\), so we see that \(N\) is in the genus of \(\bE_8 \oplus [-2m]\).

The converse holds by \cite[Prop.~1.15.1]{Nikulin:int.sym.bilinear.forms}, because
\(S_X\) is unique in its genus, \(K\) is uniquely determined and \(\OO(\bM) \rightarrow \OO(\bM^\sharp)\) is surjective (see for instance \cite[p. 388]{Barth.Peters:aut.enriques}).
\endproof

Note that a lattice \(N' \cong N(2)\), with \(N\) an even lattice, does not contain vectors of square \(-2\). 
Therefore, \autoref{lem:emb.M-->S_X.for.T_X=U+4m} essentially says that the terms in the sum of \autoref{thm:Shimada.Veniani} are in one-to-one correspondence with the lattices in the genus of \(\bE_8 \oplus [-2m]\).
In particular, one of them corresponds to \(\bE_8 \oplus [-2m]\) itself, which we presently consider more in detail.

Barth--Peters introduced a \(2\)-dimensional family of K3 surfaces, whose general element \(\fX\) has transcendental lattice \(T_\fX \cong \bU \oplus \bU(2)\) and Néron--Severi lattice \(S_\fX \cong \bU(2) \oplus \bE_8^2\). 
Ohashi~\cite[Remark~4.9(2)]{Ohashi:number.Enriques} proved that \(|{\Enr(\fX)}| = 1\).
The coinvariant lattice of an Enriques involution on \(\fX\) is isomorphic to \(\bE_8(2)\). 

In the situation of \autoref{lem:emb.M-->S_X.for.T_X=U+4m}, if \(\bE_8(2)\) embeds into \((S_X)_\varepsilon \cong \iota(\bM)^\perp\), then \(\iota(\bM)^\perp \cong \bE_8(2) \oplus [-4m]\). 
This motivates the following definition.

\begin{definition}
We say that an Enriques involution \(\varepsilon \in \Aut(X)\) on \(X\) is \emph{of Barth--Peters type} if \((S_X)_\varepsilon \cong \bE_8(2) \oplus [-4m]\).
The corresponding Enriques quotient is also called \emph{of Barth--Peters type}.
\end{definition}

The following lemma provides the number of Enriques quotients of Barth--Peters type up to isomorphisms.

\begin{lemma} \label{lem:BP.quotients}
If \(X\) is a K3 surface with transcendental lattice \(T_X \cong \bU \oplus [4m]\), \(m \geq 1\), 
and \(\iota\colon \bM \hookrightarrow S_X\) is a primitive embedding with \(\iota(\bM)^\perp \cong \bE_8(2) \oplus [-4m]\), then it holds
\[
  |{\OO_\hodge^\sharp(T_X)}\backslash{\OO(T_X^\sharp)}/{\OO^\sharp(S_X,\iota)}| = 2^{\omega-1},
\]
where \(\omega\) is the number of prime divisors of~\(2m\).
\end{lemma}

\proof 
As in the proof of \autoref{lem:BP.fibrations}, it holds \(|{\OO(T_X^\sharp)}| = 2^{\omega}\).
As \(\rank(T_X)\) is odd, it holds \(\OO_\hodge^\sharp(T_X) = \set{\pm \id}\) (see for instance \cite[Cor.~3.3.5]{Huybrechts:lectures.K3}). 
Note, moreover, that \(\id \neq -\id\) in~\(T_X^\sharp\).

We now want to determine \(\OO_\hodge^\sharp(S_X,\iota)\) using the identification \eqref{eq:L^sharp=Gamma^perp/Gamma}.
Let \(s \in \bE_8(2) \oplus [-4m]\) be the generator of the copy of \([-4m]\), \(H\) be the gluing subgroup, \(\gamma \colon H \rightarrow H'\) be the gluing isometry, and \(\Gamma\) the gluing graph of \(\iota\) (see \autoref{sec:lattices}). By the identification \eqref{eq:L^sharp=Gamma^perp/Gamma}, it holds \(|H| = |H'| = 2^9\). 
Therefore, \(S_X^\sharp \cong \Gamma^\perp/\Gamma\) is generated by an element of the form \((\alpha,s/4m)\), with \(\alpha \in \bM^\sharp\).

Recall now that \(\OO(\bM) \rightarrow \OO(\bM^\sharp)\) is surjective (see \cite[p. 388]{Barth.Peters:aut.enriques}) and that an isometry of a definite lattice preserves its decomposition in irreducible lattices up to order (see for instance \cite[Satz~27.2]{Kneser:quadratische.Formen}), so \(\OO(\bE_8(2) \oplus [-4m]) \cong \OO(\bE_8(2)) \times \OO([-4m])\).
These facts imply that the group \(\OO_\hodge(S_X,\iota)\) can only act as \(\pm \id\) on \((\alpha,s/4m)\), i.e. \(\OO_\hodge^\sharp(S_X,\iota) = \set{\pm \id}\).
Therefore, we have
\[
    |{\OO_\hodge^\sharp(T_X)}\backslash{\OO(T_X^\sharp)}/{\OO^\sharp(S_X,\iota)}| = |{\OO(T_X^\sharp)}/{\set{\pm \id}}| = |{\OO(T_X^\sharp)}| / |{\set{\pm \id}}| = 2^{\omega-1}. \qedhere
\]
\endproof

\begin{remark} \label{rmk:pullbacks}
Recall that any elliptic pencil on an Enriques surface has exactly two multiple fibers \(2F\), \(2F'\). The divisors \(F\) and \(F'\) are called \emph{half-pencils} (necessarily of type \(\I_m\) for some \(m \geq 0\)). 
An elliptic pencil on an Enriques surface is said to be \emph{special} if it has a \(2\)-section which is a smooth rational curve. 

As noted by Kond\={o} \cite[Lem.~2.6]{Kondo:Enriques.finite.aut}, the pullback of a special elliptic pencil induces a jacobian elliptic fibration \(\pi\) on the K3 surface \(X\).
Such pullbacks satisfy the following condition:
if the fibration \(\pi\) has exactly \(n_i\) fibres type of type \(\J_i\) (for \(i = 1,\ldots,r\)), where \(\J_i, \J_j\) are pairwise distinct Kodaira types if \(i \neq j\),
then at most two coefficients \(n_i\) can be odd; 
moreover, if \(n_i\) is odd, then \(\J_i = \I_{2m}\) for some \(m\geq 0\). 

The last sentence comes from the fact that one of the fibers of type \(\J_i\) is necessarily the pullback of a half-pencil.
\end{remark} 

\begin{remark} \label{rmk:Barth-Peters.ell.fibr}
Let \(n = 2m\) be an even integer and consider one of the \(2^{\omega-1}\) elliptic fibrations with two fibers of type~\(\II^*\) given in \autoref{lem:BP.fibrations}. By the construction of \autoref{rmk:Hulek.Schutt.II^*}, we obtain one of the \(2^{\omega-1}\) Enriques quotients \(Y\) of Barth--Peters type of \autoref{lem:BP.quotients}.
In the notation of \autoref{rmk:Hulek.Schutt.II^*}, the vector
\[
    R \coloneqq (m + 1, 2, -4, -5, -7, -10, -8, -6, -4, -2, -2, -3, -4, -6, -5, -4, -3, -2, 1)
\]
has square \(-2\), satisfies \(R \cdot F = 2\) and has intersection number \(1\) with \(e_3\) and \(e_{18}\). Therefore, it represents a smooth rational curve. Moreover, \(R \cdot \varepsilon(R) = 0\).

Thus, the surface \(Y\) contains ten smooth rational curves which are the images of \(R\) and of the components of the fibers of type \(\II^*\). They form the following dual graph, where the white vertex represents the image of \(R\), which is a \(4\)-section of the highlighted elliptic pencil.

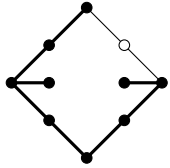
\begin{center}
\begin{tikzpicture}[scale=0.5]
    \draw[very thick] (-2, 0) -- (-1, 0);
    \draw[very thick] ( 1, 0) -- ( 2, 0);
    
    \draw[very thick] (-2, 0) -- (-1, 1);
    \draw[very thick] (-1, 1) -- ( 0, 2);
    \draw[very thick] (-2, 0) -- (-1,-1);
    \draw[very thick] (-1,-1) -- ( 0,-2);
    
    \draw             ( 2, 0) -- ( 1, 1);
    \draw             ( 1, 1) -- ( 0, 2);
    \draw[very thick] ( 2, 0) -- ( 1,-1);
    \draw[very thick] ( 1,-1) -- ( 0,-2);
    
    \draw[fill=black] (-2, 0) circle (4pt) node {};
    \draw[fill=black] (-1, 0) circle (4pt) node {};
    \draw[fill=black] (-1, 1) circle (4pt) node {};
    \draw[fill=black] ( 0, 2) circle (4pt) node {};
    \draw[fill=white] ( 1, 1) circle (4pt) node {};
    \draw[fill=black] ( 2, 0) circle (4pt) node {};
    \draw[fill=black] ( 1, 0) circle (4pt) node {};
    \draw[fill=black] (-1,-1) circle (4pt) node {};
    \draw[fill=black] ( 1,-1) circle (4pt) node {};
    \draw[fill=black] ( 0,-2) circle (4pt) node {};
\end{tikzpicture}
\end{center}

This graph appears in \cite[Thm.~1.7(i)]{Kondo:Enriques.finite.aut} and is related to the fact that \(Y\) has a cohomologically trivial automorphism (such automorphisms were studied by Mukai and Namikawa~\cite{Mukai:numerically.trivial.inv.Kummer.type,Mukai.Namikawa:aut.Enriques.which.act.trivially}).

On the above graph we can recognize three more special elliptic pencils up to symmetries (dashed lines indicates half-pencils):

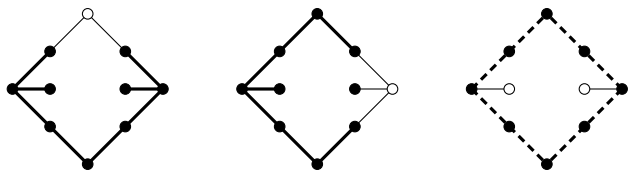
\begin{center}
\begin{tikzpicture}[scale=0.5]
    \draw[very thick] (-2, 0) -- (-1, 0);
    \draw[very thick] ( 1, 0) -- ( 2, 0);
    
    \draw[very thick] (-2, 0) -- (-1, 1);
    \draw             (-1, 1) -- ( 0, 2);
    \draw[very thick] (-2, 0) -- (-1,-1);
    \draw[very thick] (-1,-1) -- ( 0,-2);
    
    \draw[very thick] ( 2, 0) -- ( 1, 1);
    \draw             ( 1, 1) -- ( 0, 2);
    \draw[very thick] ( 2, 0) -- ( 1,-1);
    \draw[very thick] ( 1,-1) -- ( 0,-2);
    
    \draw[fill=black] (-2, 0) circle (4pt) node {};
    \draw[fill=black] (-1, 0) circle (4pt) node {};
    \draw[fill=black] (-1, 1) circle (4pt) node {};
    \draw[fill=white] ( 0, 2) circle (4pt) node {};
    \draw[fill=black] ( 1, 1) circle (4pt) node {};
    \draw[fill=black] ( 2, 0) circle (4pt) node {};
    \draw[fill=black] ( 1, 0) circle (4pt) node {};
    \draw[fill=black] (-1,-1) circle (4pt) node {};
    \draw[fill=black] ( 1,-1) circle (4pt) node {};
    \draw[fill=black] ( 0,-2) circle (4pt) node {};
\end{tikzpicture}
\qquad 
\begin{tikzpicture}[scale=0.5]
    \draw[very thick] (-2, 0) -- (-1, 0);
    \draw ( 1, 0) -- ( 2, 0);
    
    \draw[very thick] (-2, 0) -- (-1, 1);
    \draw[very thick] (-1, 1) -- ( 0, 2);
    \draw[very thick] (-2, 0) -- (-1,-1);
    \draw[very thick] (-1,-1) -- ( 0,-2);
    
    \draw ( 2, 0) -- ( 1, 1);
    \draw[very thick]             ( 1, 1) -- ( 0, 2);
    \draw  ( 2, 0) -- ( 1,-1);
    \draw[very thick] ( 1,-1) -- ( 0,-2);
    
    \draw[fill=black] (-2, 0) circle (4pt) node {};
    \draw[fill=black] (-1, 0) circle (4pt) node {};
    \draw[fill=black] (-1, 1) circle (4pt) node {};
    \draw[fill=black] ( 0, 2) circle (4pt) node {};
    \draw[fill=black] ( 1, 1) circle (4pt) node {};
    \draw[fill=white] ( 2, 0) circle (4pt) node {};
    \draw[fill=black] ( 1, 0) circle (4pt) node {};
    \draw[fill=black] (-1,-1) circle (4pt) node {};
    \draw[fill=black] ( 1,-1) circle (4pt) node {};
    \draw[fill=black] ( 0,-2) circle (4pt) node {};
\end{tikzpicture}
\qquad
\begin{tikzpicture}[scale=0.5]
    \draw (-2, 0) -- (-1, 0);
    \draw ( 1, 0) -- ( 2, 0);
    
    \draw[very thick,densely dashed] (-2, 0) -- (-1, 1);
    \draw[very thick,densely dashed] (-1, 1) -- ( 0, 2);
    \draw[very thick,densely dashed] (-2, 0) -- (-1,-1);
    \draw[very thick,densely dashed] (-1,-1) -- ( 0,-2);
    
    \draw[very thick,densely dashed] ( 2, 0) -- ( 1, 1);
    \draw[very thick,densely dashed] ( 1, 1) -- ( 0, 2);
    \draw[very thick,densely dashed] ( 2, 0) -- ( 1,-1);
    \draw[very thick,densely dashed] ( 1,-1) -- ( 0,-2);
    
    \draw[fill=black] (-2, 0) circle (4pt) node {};
    \draw[fill=white] (-1, 0) circle (4pt) node {};
    \draw[fill=black] (-1, 1) circle (4pt) node {};
    \draw[fill=black] ( 0, 2) circle (4pt) node {};
    \draw[fill=black] ( 1, 1) circle (4pt) node {};
    \draw[fill=black] ( 2, 0) circle (4pt) node {};
    \draw[fill=white] ( 1, 0) circle (4pt) node {};
    \draw[fill=black] (-1,-1) circle (4pt) node {};
    \draw[fill=black] ( 1,-1) circle (4pt) node {};
    \draw[fill=black] ( 0,-2) circle (4pt) node {};
\end{tikzpicture}
\end{center}
In our case we retrieve the jacobian elliptic fibrations on \(X\) with respectively two fibres of type~\(\I_4^*\), two fibres of type~\(\III^*\), and one fibre of type \(\I_{16}\) (hence the corresponding root lattices \(W_\rootlattice\) of the frame contain the sublattices \(\bD_8^2\), \(\bE_7^2\) and \(\bA_{15}\), respectively).
\end{remark}

\section{The three pencils}

This section is divided into three subsections, in which we study K3 surface \(X\) with transcendental lattice \(T_X \cong \bU \oplus [4]\) (\autoref{sec:Kondo-I}), \(T_X \cong \bU \oplus [8]\) (\autoref{sec:Kondo-II}), and \(T_X \cong \bU \oplus [12]\) (\autoref{sec:Apery.Fermi}).
In each case we determine \(|{\Enr(X)}|\) and \(|\cJ_X/{\Aut(X)}|\), then we focus on their Enriques quotients, especially those \emph{not} of Barth--Peters type (because those of Barth--Peters type were already considered in \autoref{sec:Barth-Peters.type}).

Moreover, we show that all jacobian elliptic fibrations satisfying the condition in~\autoref{rmk:pullbacks} are indeed pullbacks of elliptic pencils on some Enriques quotient.

\subsection{\texorpdfstring{Kond\={o}}{Kondo}'s pencil I} \label{sec:Kondo-I}
Let \(X\) be a K3 surface with transcendental lattice
\[
    T_X \cong \bU \oplus [4].
\]

\begin{theorem} \label{thm:Enr(X)_Kondo_I}
It holds \({|\Enr(X)|} = 1\).
\end{theorem}
\proof
The lattice \(\bE_8 \oplus [-2]\) is unique in its genus by the mass formula. By \autoref{lem:emb.M-->S_X.for.T_X=U+4m}, the sum in \autoref{thm:Shimada.Veniani} has only one term, which is equal to \(1\) by \autoref{lem:BP.quotients}.
\endproof 

Therefore, the surface \(X\) admits only one Enriques quotient \(X \rightarrow Y\). Necessarily, the Barth--Peters quotient of \autoref{lem:BP.quotients} coincides with Kond\={o}'s quotient~\cite{Kondo:Enriques.finite.aut} (in particular, \(Y\) has a finite automorphism group).
Indeed, the graph of nodal curves contained in \(Y\), which is pictured in \cite[Fig.~1.4]{Kondo:Enriques.finite.aut}, contains the Barth--Peters graph as a subgraph.
This Enriques quotient was also studied by Hulek and Schütt~\cite[§4.6]{Hulek.Schuett:Enriques.surf.jacobian.ell}.

For the sake of completeness, we enumerate all jacobian elliptic fibrations on \(X\) up to automorphisms (the same list is contained in an unpublished paper by Elkies and Schütt~\cite{Elkies.Schuett:K3.families.high.Picard}). 

\begin{proposition}
The frame genus \(\cW_X\) contains exactly \(9\) isomorphism classes, listed in \autoref{tab:genus_Kondo_I}, whose Gram matrices are contained in the arXiv ancillary file \verb+genus_Kondo_I.sage+.
Moreover, it holds
\[
    |\cJ_X/{\Aut(X)}| = 9.
\]
\end{proposition}
\proof 
It holds \(|\cJ_X/{\Aut(X)}| = |\cW_X|\) by \cite[Cor.~2.10]{Festi.Veniani:counting.ell.fibr.K3}. 
In order to determine \(\cW_X\), we apply the Kneser--Nishiyama method with \(T_0 = \bD_7\).
The list is complete because the mass formula holds:
\[
\sum_{i = 1}^{9} \frac{1}{|{\OO(W_i)}|} = \frac{642332179}{18881368343036559360000} = \mass(\cW_X). \qedhere 
\]
\endproof 

{\small
\begin{table}[t]
    \centering
    \caption{Lattices in the frame genus \(\cW_X\) of a K3 surface \(X\) with transcendental lattice \(T_X \cong \bU \oplus [4]\).}
    \label{tab:genus_Kondo_I}
    \begin{tabular}{llllllll}
    \toprule 
    \(W\) & \(N_\rootlattice\) & \(W_\rootlattice\) & \(W/W_\rootlattice\) & \(|\Delta(W)|\) & \(|{\OO(W)}|\) & \(|{\fr_X^{-1}(W)}|\) & Rmk. \\
    \midrule
    \(W_1\) & \(\bD_{16}\bE_8\)   & \(\bD_9\bE_8\)    & \(0\) & \(384\) & \(129448569470976000\) & \(1\) & -- \\
    \(W_2\) & \(\bD_{24}\)        & \(\bD_{17}\)      & \(0\) & \(544\) & \(46620662575398912000\)  & \(1\) & -- \\
    \midrule 
    \(W_3\) & \(\bD_{10}\bE_7^2\)  & \(\bA_3 \bE_7^2\)      & \(\IZ/2\IZ\) & \(264\) & \(809053559193600\) & \(1\) & \ref{rmk:Barth-Peters.ell.fibr} \\
    \(W_4\) & \(\bD_{12}^2\)       & \(\bD_5\bD_{12}\)     & \(\IZ/2\IZ\) & \(304\) & \(3767021862912000\) & \(1\) & -- \\
    \midrule 
    \(W_5\) & \(\bA_{11}\bD_7\bE_6\) & \(\bA_{11}\bE_6\)  & \(\IZ/3\IZ\) & \(204\) & \(49662885888000\) & \(1\) & -- \\ 
    \midrule 
    \(W_6\) & \(\bA_{15}\bD_9\)   & \(\bA_1^2\bA_{15}\)    & \(\IZ/4\IZ\) & \(244\) & \(334764638208000\) & \(1\) & \ref{rmk:Barth-Peters.ell.fibr} \\
    \midrule 
    \(W_7\) & \(\bE_8^3\)          & \(\bE_8^2\)            & \(\IZ\) & \(480\) & \(1941728542064640000\) & \(1\) & \ref{rmk:Barth-Peters.ell.fibr} \\ 
    \midrule 
    \(W_8\) & \(\bD_8^3\)          & \(\bD_8^2\)            & \(\IZ \oplus (\IZ/2\IZ)\) & \(224\) & \(106542032486400\) & \(1\) & \ref{rmk:Barth-Peters.ell.fibr} \\ 
    \(W_9\) & \(\bD_{16}\bE_8\)   & \(\bD_{16}\)          & \(\IZ \oplus (\IZ/2\IZ)\) & \(480\) & \(1371195958099968000\) & \(1\) & -- \\ 
    \bottomrule
    \end{tabular}
\end{table}
}

\subsection{\texorpdfstring{Kond\={o}}{Kondo}'s pencil II} \label{sec:Kondo-II}

Let \(X\) be a K3 surface with transcendental lattice
\[
    T_X \cong \bU \oplus [8].
\]

\begin{theorem} \label{thm:Enr(X)_Kondo_II}
It holds \({|\Enr(X)|} = 2\).
\end{theorem}
\proof
There is only one more lattice in the genus of \(\bE_8 \oplus [-4]\), namely \(\bD_9\), as the mass formula shows.
By \autoref{lem:emb.M-->S_X.for.T_X=U+4m} the sum in \autoref{thm:Shimada.Veniani} has two terms, both equal to \(1\) as it holds \(\OO_\hodge^\sharp(T_X) = \OO(T_X^\sharp)\).
\endproof 

By \autoref{lem:BP.quotients}, one of the two Enriques quotient, say \(X \rightarrow Y'\), is of Barth--Peters type. The corresponding coinvariant lattice is isometric to \(\bE_8(2) \oplus [-8]\) and contains \(240\) vectors of square \(-4\).

By Kond\={o}'s classification, the surface \(X\) admits an Enriques quotient \(X \rightarrow Y''\) with finite automorphism group. 
Kond\={o}'s quotient \(Y''\) was also studied by Hulek and Schütt \cite[§4.7 and §4.8]{Hulek.Schuett:Enriques.surf.jacobian.ell}.
We argue that \(Y'\) is not isomorphic to \(Y''\).

Geometrically, this follows from the fact that \(Y''\) contains exactly \(12\) smooth rational curves whose dual graph is pictured on \cite[p.~207, Fig.~2.4]{Kondo:Enriques.finite.aut}. This dual graph does not contain the graph pictured in \autoref{rmk:Barth-Peters.ell.fibr} as a subgraph.

Algebraically, we can distinguish the two quotients in the following way.
\begin{enumerate}
    \item The surface \(X\) contains \(24\) smooth rational curves \(F_1^+,F_1^-,\ldots,F_{12}^+,F_{12}^-\), which intersect as in \cite[p.~207, Fig.~2.3]{Kondo:Enriques.finite.aut} and generate the Néron--Severi lattice \(S_X\).
    \item Kond\={o}'s Enriques involution exchanges \(F_i^+\) with \(F_i^-\), \(i = 1,\ldots,12\).
    \item Computing explicitly the coinvariant lattice of Kond\={o}'s Enriques involution in \(S_X\), we see that it contains \(144\) vectors of square \(-4\), so it must be isomorphic to \(\bD_9(2)\). In particular, Kond\={o}'s quotient is not of Barth--Peters type.
\end{enumerate}

We now enumerate all jacobian elliptic fibrations on~\(X\) up to automorphisms.

{\small
\begin{table}[t]
    \centering
    \caption{Lattices in the frame genus \(\cW_X\) of a K3 surface \(X\) with transcendental lattice \(T_X \cong \bU \oplus [8]\).}
    \label{tab:genus_Kondo_II}
    \begin{tabular}{llllllll}
    \toprule 
    \(W\) & \(N_\rootlattice\) & \(W_\rootlattice\) & \(W/W_\rootlattice\) & \(|\Delta(W)|\) & \(|{\OO(W)}|\) & \(|{\fr_X^{-1}(W)}|\) & Rmk. \\
    \midrule
    \(W_1\) & \(\bA_7^2\bD_5^2\) & \(\bA_7\bD_5^2\) & \(\IZ/4\IZ\) & \(136\) & \(594542592000\) & \(1\) & \ref{rmk:Kondo-II.ell.fibr} \\
    \midrule
    \(W_2\) & \(\bA_{11}\bD_7\bE_6\) & \(\bA_3\bD_7\bE_6\) & \(\IZ\) & \(168\) & \(802632499200\) & \(1\) & -- \\
    \(W_3\) & \(\bA_{12}^2\) & \(\bA_{12}\bA_4\) & \(\IZ \) & \(176\) & \(1494484992000\) & \(1\) & -- \\
    \(W_4\) & \(\bA_{15}\bD_9\) & \(\bA_7\bD_9\) & \(\IZ\) & \(200\) & \(7491236659200\) & \(1\) & -- \\
    \(W_5\) & \(\bA_{17}\bE_7\) & \(\bA_9\bE_7\) & \(\IZ\) & \(216\) & \(21069103104000\) & \(1\) & -- \\
    \(W_6\) & \(\bA_{24}\) & \(\bA_{16}\) & \(\IZ\) & \(272\) & \(711374856192000\) & \(1\) & -- \\
    \(W_7\) & \(\bD_{16}\bE_8\) & \(\bD_8\bE_8\) & \(\IZ\) & \(342\) & \(7191587192832000\) & \(1\) & -- \\
    \(W_8\) & \(\bD_{24}\) & \(\bD_{16}\) & \(\IZ\) & \(480\) & \(1371195958099968000\) & \(1\) & -- \\
    \(W_9\) & \(\bE_8^3\) & \(\bE_8^2\) & \(\IZ\) & \(480\) & \(1941728542064640000\) & \(1\) & \ref{rmk:Barth-Peters.ell.fibr} \\
    \midrule
    \(W_{10}\) & \(\bA_9^2 \bD_6\) & \(\bA_1\bA_9\bD_6\) & \(\IZ \oplus (\IZ/2\IZ)\) & \(152\) & \(334430208000\) & \(1\) & -- \\ 
    \(W_{11}\) & \(\bD_8^3\) & \(\bD_8^2\) & \(\IZ \oplus (\IZ/2\IZ)\) & \(224\) & \(106542032486400\) & \(1\) & \ref{rmk:Kondo-II.ell.fibr}, \ref{rmk:W11,W12} \\
    \(W_{12}\) & \(\bD_8^3\) & \(\bD_8^2\) & \(\IZ \oplus (\IZ/2\IZ)\) & \(224\) & \(106542032486400\) & \(1\) & \ref{rmk:Barth-Peters.ell.fibr}, \ref{rmk:W11,W12} \\
    \(W_{13}\) & \(\bD_{10}\bE_7^2\) & \(\bA_1^2 \bE_7^2\) & \(\IZ \oplus (\IZ/2\IZ)\) & \(256\) & \(134842259865600\) & \(1\) & \ref{rmk:Barth-Peters.ell.fibr} \\
    \(W_{14}\) & \(\bD_{12}^2\) & \(\bD_{12}\bD_4\) & \(\IZ \oplus (\IZ/2\IZ)\) & \(288\) & \(376702186291200\) & \(1\) & -- \\
    \(W_{15}\) & \(\bD_{16}\bE_8\) & \(\bD_{16}\) & \(\IZ \oplus (\IZ/2\IZ)\) & \(480\) & \(1371195958099968000\) & \(1\) & -- \\
    \midrule
    \(W_{16}\) & \(\bA_8^3\) & \(\bA_8^2\) & \(\IZ \oplus (\IZ/3\IZ)\) & \(144\) & \(526727577600\) & \(1\) & \ref{rmk:Kondo-II.ell.fibr} \\
    \midrule    
    \(W_{17}\) & \(\bA_{15}\bD_9\) & \(\bA_{15}\) & \(\IZ^2 \oplus (\IZ/2\IZ)\) & \(240\) & \(83691159552000\) & \(1\) & \ref{rmk:Barth-Peters.ell.fibr} \\
    \bottomrule
    \end{tabular}
\end{table}
}

\begin{proposition}
The frame genus \(\cW_X\) contains exactly \(17\) isomorphism classes, listed in \autoref{tab:genus_Kondo_I}, whose Gram matrices are contained in the arXiv ancillary file \verb+genus_Kondo_II.sage+.
Moreover, it holds
\[
    |\cJ_X/{\Aut(X)}| = 17.
\]
\end{proposition}
\proof 
It holds \(|\cJ_X/{\Aut(X)}| = |\cW_X|\) by \cite[Cor.~2.10]{Festi.Veniani:counting.ell.fibr.K3}. 
In order to determine \(\cW_X\), we apply the Kneser--Nishiyama method with \(T_0 = \bA_7\). 
Note that there are two different primitive embeddings \(\bA_7 \hookrightarrow \bD_8\) (cf. \cite[Lem.~4.2]{Nishiyama:Jacobian.fibrations.K3.MW}), leading to two distinct frames \(W\) with \(W_\rootlattice \hookrightarrow \bD_8^3\), namely \(W_{11}\) and~\(W_{12}\) (cf. \autoref{rmk:W11,W12}).
The list is complete because the mass formula holds:
\[
\sum_{i = 1}^{17} \frac{1}{|{\OO(W_i)}|} = \frac{642332179}{73755345089986560000} = \mass(\cW_X). \qedhere 
\]
\endproof 
\begin{remark} \label{rmk:Kondo-II.ell.fibr}
The surface \(Y''\) contains \(12\) curves on whose dual graph one can recognize the following elliptic pencils (dashed lines indicate half-pencils):
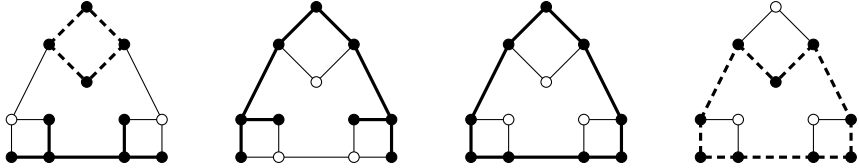
\begin{center}
\begin{tikzpicture}[scale=0.5]
    \draw (-2, 0) -- (-1, 2);
    \draw[very thick,densely dashed] (-1, 2) -- ( 0, 3);
    \draw[very thick,densely dashed] ( 1, 2) -- ( 0, 3);
    \draw ( 1, 2) -- ( 2, 0);
    \draw ( 2, 0) -- ( 2,-1);
    \draw[very thick] ( 2,-1) -- ( 1,-1);
    \draw[very thick] ( 1,-1) -- (-1,-1);
    \draw[very thick] (-1,-1) -- (-2,-1);
    \draw (-2,-1) -- (-2, 0);
    
    \draw (-2, 0) -- (-1, 0);
    \draw[very thick] (-1,-1) -- (-1, 0);
    
    \draw[very thick,densely dashed] (-1, 2) -- ( 0, 1);
    \draw[very thick,densely dashed] ( 1, 2) -- ( 0, 1);
    
    \draw ( 2, 0) -- ( 1, 0);
    \draw[very thick] ( 1,-1) -- ( 1, 0);
    
    \draw[fill=white] (-2, 0) circle (4pt) node {};
    \draw[fill=black] (-1, 0) circle (4pt) node {};
    \draw[fill=black] ( 1, 0) circle (4pt) node {};
    \draw[fill=white] ( 2, 0) circle (4pt) node {};
    \draw[fill=black] (-2,-1) circle (4pt) node {};
    \draw[fill=black] (-1,-1) circle (4pt) node {};
    \draw[fill=black] ( 1,-1) circle (4pt) node {};
    \draw[fill=black] ( 2,-1) circle (4pt) node {};

    \draw[fill=black] (-1, 2) circle (4pt) node {};
    \draw[fill=black] ( 0, 1) circle (4pt) node {};
    \draw[fill=black] ( 0, 3) circle (4pt) node {};
    \draw[fill=black] ( 1, 2) circle (4pt) node {};
\end{tikzpicture}
\qquad
\begin{tikzpicture}[scale=0.5]
    \draw[very thick] (-2, 0) -- (-1, 2);
    \draw[very thick] (-1, 2) -- ( 0, 3);
    \draw[very thick] ( 0, 3) -- ( 1, 2);
    \draw[very thick] ( 1, 2) -- ( 2, 0);
    \draw[very thick] ( 2, 0) -- ( 2,-1);
    \draw ( 2,-1) -- ( 1,-1);
    \draw ( 1,-1) -- (-1,-1);
    \draw (-1,-1) -- (-2,-1);
    \draw[very thick] (-2,-1) -- (-2, 0);
    
    \draw[very thick] (-2, 0) -- (-1, 0);
    \draw (-1,-1) -- (-1, 0);
    
    \draw (-1, 2) -- ( 0, 1);
    \draw ( 1, 2) -- ( 0, 1);
    
    \draw[very thick] ( 2, 0) -- ( 1, 0);
    \draw ( 1,-1) -- ( 1, 0);
    
    \draw[fill=black] (-2, 0) circle (4pt) node {};
    \draw[fill=black] (-1, 0) circle (4pt) node {};
    \draw[fill=black] ( 1, 0) circle (4pt) node {};
    \draw[fill=black] ( 2, 0) circle (4pt) node {};
    \draw[fill=black] (-2,-1) circle (4pt) node {};
    \draw[fill=white] (-1,-1) circle (4pt) node {};
    \draw[fill=white] ( 1,-1) circle (4pt) node {};
    \draw[fill=black] ( 2,-1) circle (4pt) node {};

    \draw[fill=black] (-1, 2) circle (4pt) node {};
    \draw[fill=white] ( 0, 1) circle (4pt) node {};
    \draw[fill=black] ( 0, 3) circle (4pt) node {};
    \draw[fill=black] ( 1, 2) circle (4pt) node {};
\end{tikzpicture}
\qquad
\begin{tikzpicture}[scale=0.5]
    \draw[very thick] (-2, 0) -- (-1, 2);
    \draw[very thick] (-1, 2) -- ( 0, 3);
    \draw[very thick] ( 0, 3) -- ( 1, 2);
    \draw[very thick] ( 1, 2) -- ( 2, 0);
    \draw[very thick] ( 2, 0) -- ( 2,-1);
    \draw[very thick] ( 2,-1) -- ( 1,-1);
    \draw[very thick] ( 1,-1) -- (-1,-1);
    \draw[very thick] (-1,-1) -- (-2,-1);
    \draw[very thick] (-2,-1) -- (-2, 0);
    
    \draw (-2, 0) -- (-1, 0);
    \draw (-1,-1) -- (-1, 0);
    
    \draw (-1, 2) -- ( 0, 1);
    \draw ( 1, 2) -- ( 0, 1);
    
    \draw ( 2, 0) -- ( 1, 0);
    \draw ( 1,-1) -- ( 1, 0);
    
    \draw[fill=black] (-2, 0) circle (4pt) node {};
    \draw[fill=white] (-1, 0) circle (4pt) node {};
    \draw[fill=white] ( 1, 0) circle (4pt) node {};
    \draw[fill=black] ( 2, 0) circle (4pt) node {};
    \draw[fill=black] (-2,-1) circle (4pt) node {};
    \draw[fill=black] (-1,-1) circle (4pt) node {};
    \draw[fill=black] ( 1,-1) circle (4pt) node {};
    \draw[fill=black] ( 2,-1) circle (4pt) node {};
    
    \draw[fill=black] (-1, 2) circle (4pt) node {};
    \draw[fill=white] ( 0, 1) circle (4pt) node {};
    \draw[fill=black] ( 0, 3) circle (4pt) node {};
    \draw[fill=black] ( 1, 2) circle (4pt) node {};
\end{tikzpicture}
\qquad
\begin{tikzpicture}[scale=0.5]
    \draw[very thick,densely dashed] (-2, 0) -- (-1, 2);
    \draw (-1, 2) -- ( 0, 3);
    \draw ( 0, 3) -- ( 1, 2);
    \draw[very thick,densely dashed] ( 1, 2) -- ( 2, 0);
    \draw[very thick,densely dashed] ( 2, 0) -- ( 2,-1);
    \draw[very thick,densely dashed] ( 2,-1) -- ( 1,-1);
    \draw[very thick,densely dashed] ( 1,-1) -- (-1,-1);
    \draw[very thick,densely dashed] (-1,-1) -- (-2,-1);
    \draw[very thick,densely dashed] (-2, 0) -- (-2,-1);
    
    \draw (-2, 0) -- (-1, 0);
    \draw (-1,-1) -- (-1, 0);
    
    \draw[very thick,densely dashed] (-1, 2) -- ( 0, 1);
    \draw[very thick,densely dashed] ( 1, 2) -- ( 0, 1);
    
    \draw ( 2, 0) -- ( 1, 0);
    \draw ( 1,-1) -- ( 1, 0);
    
    \draw[fill=black] (-2, 0) circle (4pt) node {};
    \draw[fill=white] (-1, 0) circle (4pt) node {};
    \draw[fill=white] ( 1, 0) circle (4pt) node {};
    \draw[fill=black] ( 2, 0) circle (4pt) node {};
    \draw[fill=black] (-2,-1) circle (4pt) node {};
    \draw[fill=black] (-1,-1) circle (4pt) node {};
    \draw[fill=black] ( 1,-1) circle (4pt) node {};
    \draw[fill=black] ( 2,-1) circle (4pt) node {};
    
    \draw[fill=black] (-1, 2) circle (4pt) node {};
    \draw[fill=black] ( 0, 1) circle (4pt) node {};
    \draw[fill=white] ( 0, 3) circle (4pt) node {};
    \draw[fill=black] ( 1, 2) circle (4pt) node {};
\end{tikzpicture}
\end{center}

The first three pencils are special pencils and correspond to the elliptic fibrations on \(X\) with frames \(W_1\), \(W_{11}\) (see \autoref{rmk:W11,W12}) and \(W_{16}\), respectively. 

The fourth pencil is not special: the highlighted curves on \(Y''\) form a half-pencil and the white vertices represent \(4\)-sections. 
Indeed, the pullback on \(X\) correspond to an elliptic fibration with a fiber of type \(\I_{18}\), namely
\[
    F_1^+ + F_4^- + F_3^- + F_5^- + F_6^- + F_7^- + F_9^- + F_{10}^- + F_{11}^- + F_1^- + F_4^+ + F_3^+ + F_5^+ + F_6^+ + F_7^+ + F_9^+ + F_{10}^+ + F_{11}^+.
\]
This fibration is not jacobian, as it does not appear in \autoref{tab:genus_Kondo_II}.
\end{remark}

\begin{remark} \label{rmk:W11,W12}
The two frames \(W_{11}\) and \(W_{12}\) are not isometric, but they can be distinguished neither by the pair \((W_\rootlattice,W/W_\rootlattice)\) nor by their number of automorphisms \(|{\OO(W)}|\). 

Using the command \verb+is_globally_equivalent_to+ of the Sage class \verb+QuadraticForm+,
we can check that the frame \(W_{11}\) corresponds to the fibration with fiber
\begin{align*}
    F & = F_6^+ + F_8^- + 2F_5^+ + 2F_3^+ + 2F_2^+ + 2F_1^+ + 2F_{11}^+ + F_{10}^+ + F_{12}^+ \\
     & = F_6^- + F_8^+ + 2F_5^- + 2F_3^- + 2F_2^- + 2F_1^- + 2F_{11}^- + F_{10}^- + F_{12}^-.
\end{align*}
which is the pullback of the second special pencil on \(Y''\) listed in \autoref{rmk:Kondo-II.ell.fibr}.

On the other hand, with the same command we can check that the frame \(W_{12}\) corresponds to the fibration with fiber
\[
    F_2^+ + F_5^+ + 2F_3^+ + 2F_4^+ + 2F_1^- + 2F_2^- + 2F_3^- + F_4^- + F_5^-,
\]
which is then the pullback of a special pencil on \(Y'\) (cf.~\autoref{rmk:Barth-Peters.ell.fibr}).
\end{remark}

\subsection{Apéry--Fermi pencil} \label{sec:Apery.Fermi}

Let \(X\) be a K3 surface with transcendental lattice
\[
    T_X \cong \bU \oplus [12].
\]

The classification of the jacobian elliptic fibrations on \(X\) was carried out by Bertin and Lecacheux~\cite{Bertin.Lecacheux:apery-fermi.2-isogenies} and then refined in~\cite{Festi.Veniani:counting.ell.fibr.K3}.
For the reader's convenience we reproduce in \autoref{tab:genus_Apery_Fermi} the same table as \cite[Table~7]{Festi.Veniani:counting.ell.fibr.K3}.

{\small
\begin{table}[t]
    \centering
    \caption{Lattices in the frame genus \(\cW_X\) of a K3 surface \(X\) with transcendental lattice \(T_X \cong \bU \oplus [12]\), numbered according to Bertin and Lecacheux (cf. \cite[Tables~2 and 3]{Bertin.Lecacheux:apery-fermi.2-isogenies}).}
    \label{tab:genus_Apery_Fermi}
    \begin{tabular}{llllllll}
    \toprule 
    \(W\) & \(N_\rootlattice\) & \(W_\rootlattice\) & \(W/W_\rootlattice\) & \(|\Delta(W)|\) & \(|{\OO(W)}|\) & \(|{\fr_X^{-1}(W)}|\) & Rmk. \\
    \midrule
    \(W_3\) & \(\bD_{16}\bE_8\) & \(\bD_{11}\bE_6\)     & \(0\)     & \(292\) & \(8475799191552000\) & \(1\) & -- \\
    \(W_1\) & \(\bE_8^3\) & \(\bA_3\bE_6\bE_8\)         & \(0\)     & \(324\) & \(3467372396544000\) & \(1\) & -- \\
    \midrule
    \(W_7\) & \(\bD_{10}\bE_7^2\) & \(\bA_5\bD_5\bE_7\) & \(\IZ/2\IZ\) & \(196\) & \(16052649984000\) & \(1\) & -- \\
    \midrule
    \(W_{20}\) & \(\bA_{11}\bD_7\bE_6\) & \(\bA_1^2\bA_2^2\bA_{11}\) & \(\IZ/6\IZ\) & \(148\)
& \(551809843200\) & \(1\) & \ref{rmk:Apery-Fermi.ell.fibr} \\
    \midrule
    \(W_{27}\) & \(\bA_7^2\bD_5^2\) & \(\bA_4\bA_7\bD_5\) & \(\IZ\) & \(116\) & \(18579456000\) & \(2\) & -- \\
    \(W_{21}\) & \(\bA_{11}\bD_7\bE_6\) & \(\bA_1^2\bA_8\bE_6\) & \(\IZ\) & \(148\) & \(300987187200\) & \(1\) & -- \\
    \(W_{18}\) & \(\bA_{15}\bD_9\) & \(\bA_{12}\bD_4\)  & \(\IZ\) & \(180\) & \(4782351974400\) & \(1\) & -- \\
    \(W_{13}\) & \(\bD_{12}^2\) & \(\bD_9\bD_7\)        & \(\IZ\)   & \(228\) & \(119859786547200\) & \(1\) & -- \\
    \(W_5\) & \(\bD_{16}\bE_8\) & \(\bA_3 \bD_{13}\)    & \(\IZ\)   & \(324\) & \(2448564210892800\) & \(1\) & -- \\
    \(W_6\) & \(\bD_{16}\bE_8\) & \(\bD_8 \bE_8\)       & \(\IZ\)   & \(352\) & \(14383174385664000\) & \(1\) & -- \\
    \(W_2\) & \(\bE_8^3\) & \(\bE_8^2\)                 & \(\IZ\)   & \(480\) & \(1941728542064640000\) & \(2\) & \ref{rmk:Barth-Peters.ell.fibr} \\  % questa è E8+E8+[-12]
    \(W_{12}\) & \(\bD_{24}\) & \(\bD_{16}\)            & \(\IZ\)   & \(480\) & \(2742391916199936000\) & \(1\) & -- \\
    \midrule 
    \(W_{15}\) & \(\bD_8^3\) & \(\bA_3 \bD_5 \bD_8\)    & \(\IZ \oplus (\IZ/2\IZ)\) & \(164\) & \(951268147200\) & \(1\) & -- \\
    \(W_8\) & \(\bD_{10}\bE_7^2\) & \(\bA_1\bA_5\bD_{10}\) & \(\IZ \oplus (\IZ/2\IZ)\) & \(212\) & \(10701766656000\) & \(1\) & -- \\
    \(W_{16}\) & \(\bD_8^3\) & \(\bD_8^2\)              & \(\IZ \oplus (\IZ/2\IZ)\) & \(224\) & \(106542032486400\) & \(2\) & \ref{rmk:Barth-Peters.ell.fibr} \\
    \(W_9\) & \(\bD_{10}\bE_7^2\) & \(\bA_1^2\bE_7^2\) & \(\IZ \oplus (\IZ/2\IZ)\) & \(256\) & \(269684519731200\) & \(1\) & \ref{rmk:Barth-Peters.ell.fibr} \\
    \(W_{14}\) & \(\bD_{12}^2\) & \(\bD_4 \bD_{12}\)    & \(\IZ \oplus (\IZ/2\IZ)\) & \(288\) & \(753404372582400\) & \(1\) & -- \\
    \(W_4\) & \(\bD_{16}\bE_8\) & \(\bD_{16}\)          & \(\IZ \oplus (\IZ/2\IZ)\) & \(480\) & \(1371195958099968000\) & \(2\) & -- \\
    \midrule
    \(W_{19}\) & \(\bE_6^4\) & \(\bA_2^2\bE_6^2\)       & \(\IZ \oplus (\IZ/3\IZ)\) & \(156\) & \(773967052800\) & \(1\) & \ref{rmk:Apery-Fermi.ell.fibr} \\
    \midrule
    \(W_{26}\) & \(\bA_7^2\bD_5^2\) & \(\bA_1^2\bA_7^2\) & \(\IZ \oplus (\IZ/4\IZ)\) & \(116\) & \(52022476800\) & \(1\) & \ref{rmk:Apery-Fermi.ell.fibr} \\
    \midrule
    \(W_{25}\) & \(\bA_9^2\bD_6\) & \(\bA_6\bA_9\)      & \(\IZ^2\) & \(132\) & \(73156608000\) & \(1\) & -- \\
    \(W_{22}\) & \(\bA_{11}\bD_7\bE_6\) & \(\bA_8\bD_7\) & \(\IZ^2\) & \(156\) & \(234101145600\) & \(2\) & -- \\
    \(W_{10}\) & \(\bD_{10}\bE_7^2\) & \(\bA_1\bD_7\bE_7\) & \(\IZ^2\) & \(212\) & \(7491236659200\) & \(1\) & -- \\
    \(W_{11}\) & \(\bA_{17}\bE_7\) & \(\bA_1\bA_{14}\)  & \(\IZ^2\) & \(212\) & \(10461394944000\) & \(1\) & -- \\
    \midrule 
    \(W_{24}\) & \(\bD_6^4\) & \(\bA_3\bD_6^2\)         & \(\IZ^2 \oplus (\IZ/2\IZ)\) & \(132\) & \(101921587200\) & \(1\) & \ref{rmk:Apery-Fermi.ell.fibr} \\
    \(W_{23}\) & \(\bA_{11}\bD_7\bE_6\) & \(\bA_{11}\bD_4\) & \(\IZ^2 \oplus (\IZ/2\IZ)\) & \(156\) &
\(367873228800\) & \(1\) & -- \\
    \(W_{17}\) & \(\bA_{15}\bD_9\) & \(\bA_{15}\)       & \(\IZ^2 \oplus (\IZ/2\IZ)\) & \(240\) & \(167382319104000\) & \(1\) & \ref{rmk:Barth-Peters.ell.fibr} \\
    \bottomrule
    \end{tabular}
\end{table}
}

\begin{theorem} \label{thm:Enr(X)_Apery_Fermi}
It holds \({|\Enr(X)|} = 3\).
\end{theorem}
\proof
The genus of \(\bE_8 \oplus [-6]\) contains two lattices, namely \(\bA_2 \oplus \bE_7\) and \(\bE_8 \oplus [-6]\) itself, as the mass formula shows. 
Thus, by \autoref{lem:emb.M-->S_X.for.T_X=U+4m}, the sum in \autoref{thm:Shimada.Veniani} consists of two terms, one of which is equal to~\(2\) by \autoref{lem:BP.quotients}. 

Fix a primitive embedding \(\iota \colon \bM \hookrightarrow S_X\) with \(\iota(\bM)^\perp \cong \bA_2(2) \oplus \bE_7(2)\). Note that it holds
\[
    S_X^\sharp \cong \bw^5_{2,2}\oplus \bw^{-1}_{3,1}, \quad \bA_2(2)^\sharp \cong \bv_1 \oplus \bw^{-1}_{3,1}, \quad \bE_7(2)^\sharp = 3\bu_1 \oplus \bw^1_{2,2}.
\]
Let \(H \subset \bM^\sharp\) be the gluing subgroup (see \autoref{sec:lattices}).
By the identification~\eqref{eq:L^sharp=Gamma^perp/Gamma} we have \(|H| = 2^9\).
Thus, the image \(H' \coloneqq \gamma(H) \subset (\iota(\bM)^\perp(-1))^\sharp\) of the gluing isometry is the sum of the copy of \(\bv_1\) in \(\bA_2(2)\) and the whole group \(\bE_7(2)^\sharp\) (with inverted sign).

Consider the isometry \(\alpha \in \OO(\iota(\bM)^\perp)\) defined as \(-\id\) on the copy of \(\bA_2(2)\) and as \(\id\) on the copy of \(\bE_7(2)\). 
Since the natural homomorphism \(\OO(\bM) \rightarrow \OO(\bM^\sharp)\) is surjective, \(\alpha\) extends to an isometry \(\tilde \alpha \in \OO(S_X,\iota)\) by \cite[Cor.~1.5.2]{Nikulin:int.sym.bilinear.forms}. 

By construction of \(\alpha\) and by the above description of the gluing isometry \(\gamma\), the element \(\tilde \alpha^\sharp\) acts as \(-\id\) on the \(3\)-part of \(S_X^\sharp\) and as \(\id\) on the \(2\)-part of \(S_X^\sharp\). In particular, \(\OO^\sharp(S_X,\iota)\) contains at least three different elements, namely \(\id\), \(-\id\) and \(\tilde \alpha^\sharp\).
On the other hand, \(\OO(T_X^\sharp)\) contains exactly four elements, as it is generated by multiplication by \(-1\) and by \(5\).
Therefore, we have \(\OO^\sharp(S_X,\iota) = \OO(T_X^\sharp)\), which implies \[|{\OO_\hodge^\sharp(T_X)} \backslash{\OO(T_X^\sharp)}/ {\OO^\sharp(S_X,\iota)}| = 1.\]
In total we get \(|{\Enr(X)}| = 3\).
\endproof

Let \(Y',Y'',Y'''\) be the three Enriques quotients of \(X\) up to automorphisms. 
We can suppose  that \(Y', Y''\) are of Barth--Peters type (see \autoref{sec:Barth-Peters.type}). 
Here we are interested in studying \(Y \coloneqq Y'''\). 

Peters and Stienstra \cite{Peters.Stienstra:pencil.Apery} showed that \(X\) contains \(20 + 12\) smooth rational curves, called \emph{\(L\)-lines} and \emph{\(M\)-lines}, forming a particular configuration, which we call the \emph{Peters--Stientra cube}. 
The dual graph of the \(L\)-lines is pictured in \cite[Fig.~1]{Peters.Stienstra:pencil.Apery}. 
We do not reproduce it here, but we follow the same notation. 
The intersection numbers of the \(M\)-lines are described in \cite[Lem.~1]{Peters.Stienstra:pencil.Apery}.

In order to make a connection with the construction of \autoref{rmk:Hulek.Schutt.II^*}, we first look for a fibration with two fibers of type \(\II^*\) or, equivalently, with frame \(W_{2}\) in \autoref{tab:genus_Apery_Fermi}.
We can suppose
\begin{align*} 
    F & = 2 L_{++0} + 3 M_{2-+} + 4 L_{+++} + 6 L_{+0+} + 5 L_{+-+} + 4 L_{0-+} + 3 M_{1--} + 2 L_{0+-} + L_{-+-} \\
      & = 2 L_{--0} + 3 M_{1++} + 4 L_{---} + 6 L_{0--} + 5 L_{+--} + 4 L_{+0-} + 3 M_{2+-} + 2 L_{-0+} + L_{-++},
\end{align*}
as pictured below in the Peters--Stienstra cube (note that \(M_{1++}\) and \(M_{2+-}\) and the other \(M\)-lines are not displayed):
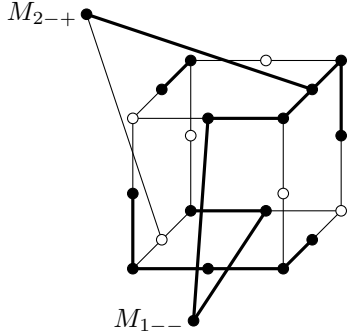
\[
\begin{tikzpicture}[x={(1cm,0cm)},y={(3.85mm, 3.85mm)},z={(0cm,1cm)},baseline=0]
    \coordinate (M2mp) at (-2, 0, 2);
    
    \draw (-1,-1,-1) -- (-1,-1, 1);
    \draw (-1,-1,-1) -- (-1, 1,-1);
    \draw[very thick] (-1,-1,-1) -- ( 1,-1,-1);
    \draw (-1,-1, 1) -- (-1, 1, 1);
    \draw (-1,-1, 1) -- ( 1,-1, 1);
    \draw (-1, 1,-1) -- (-1, 1, 1);
    \draw (-1, 1,-1) -- ( 1, 1,-1);
    \draw ( 1,-1,-1) -- ( 1,-1, 1);
    \draw ( 1,-1,-1) -- ( 1, 1,-1);
    \draw (-1, 1, 1) -- ( 1, 1, 1);
    \draw[very thick] ( 1,-1, 1) -- ( 1, 1, 1);
    \draw ( 1, 1,-1) -- ( 1, 1, 1);
    
    \draw[very thick] (-1, 1,-1) -- ( 0, 1,-1);
    \draw[very thick] ( 1, 1, 0) -- ( 1, 1, 1);
    \draw[very thick] ( 1,-1,-1) -- ( 1, 0,-1);
    \draw[very thick] (-1,-1,-1) -- (-1,-1, 0);
    \draw[very thick] ( 0,-1, 1) -- ( 1,-1, 1);
    \draw[very thick] (-1, 1, 1) -- (-1, 0, 1);
    
    \draw[very thick] ( 1, 0, 1) -- (M2mp);
    \draw (-1, 0,-1) -- (M2mp);
    
    \draw[very thick] ( 0,-1, 1) -- ( 0,-1.5,-1.5);
    \draw[very thick] ( 0, 1,-1) -- ( 0,-1.5,-1.5);
    
    \draw[fill=black] (-1,-1,-1) circle (2pt) node[above right] {};
    \draw[fill=white] (-1,-1, 1) circle (2pt) node[below] {};
    \draw[fill=black] (-1, 1,-1) circle (2pt) node[above] {};
    \draw[fill=black] (-1, 1, 1) circle (2pt) node[left] {};
    \draw[fill=black] ( 1,-1,-1) circle (2pt) node[right] {};
    \draw[fill=black] ( 1,-1, 1) circle (2pt) node[below] {};
    \draw[fill=white] ( 1, 1,-1) circle (2pt) node[above] {};
    \draw[fill=black] ( 1, 1, 1) circle (2pt) node[below right] {};
    
    \draw[fill=black] ( 0,-1,-1) circle (2pt) node[below] {};
    \draw[fill=black] ( 0,-1, 1) circle (2pt) node[below] {};
    \draw[fill=black] ( 0, 1,-1) circle (2pt) node[above] {};
    \draw[fill=white] ( 0, 1, 1) circle (2pt) node[above] {};
    
    \draw[fill=white] (-1, 0,-1) circle (2pt) node[left] {};
    \draw[fill=black] (-1, 0, 1) circle (2pt) node[left] {};
    \draw[fill=black] ( 1, 0,-1) circle (2pt) node[right] {};
    \draw[fill=black] ( 1, 0, 1) circle (2pt) node[right] {};
    
    \draw[fill=black] (-1,-1, 0) circle (2pt) node[below right] {};
    \draw[fill=white] (-1, 1, 0) circle (2pt) node[above left] {};
    \draw[fill=white] ( 1,-1, 0) circle (2pt) node[below right] {};
    \draw[fill=black] ( 1, 1, 0) circle (2pt) node[above left] {};

    \draw[fill=black] (M2mp) circle (2pt) node[left] {\(M_{2-+}\)};
    
    \draw[fill=black] ( 0,-1.5,-1.5) circle (2pt) node[left] {\(M_{1--}\)};

\end{tikzpicture} 
\]

In the coordinate system of \autoref{rmk:Hulek.Schutt.II^*}, up to substituting \(P\) with \(Q\), we can suppose that
\begin{align*}
    L_{-+0} & = O, \\
    L_{++0} & = (0,0,1,0,0,0,0,0,0,0,0,0,0,0,0,0,0,0,0), \\
    L_{--0} & = (0,0,0,0,0,0,0,0,0,0,1,0,0,0,0,0,0,0,0), \\
    L_{++-} & = (4, 4, -6, -8, -11, -16, -13, -10, -7, -4, -6, -9, -12, -18, -15, -12, -8, -4, 1).
\end{align*}
The coordinates of all other \(L\)-lines and \(M\)-lines are determined by these choices.

In order to construct the Enriques involution corresponding to \(Y\), we now consider a fibration with frame \(W_{19}\) in \autoref{tab:genus_Apery_Fermi}. 
As \((W_{19})_\rootlattice \cong \bA_2^2 \bE_6^2\), the fibration has two fibers of type \(\I_3\) (or \(\IV\)) and two fibers of type \(\IV^*\).

As pictured below in the Peters--Stienstra cube (omitting the \(M\)-lines) we choose
\begin{align*}
    F_{19} & \coloneqq L_{+-+} + L_{++-} + 2L_{+0+} + 2L_{++0} + 3L_{+++} + 2L_{0++} + L_{-++} = M_{3+-} + M_{1+-} + M_{2+-} \\
    & = L_{+--} + L_{-+-} + 2L_{0--} + 2L_{-0-} + 3L_{---} + 2L_{--0} + L_{--+} = M_{3--} + M_{1--} + M_{2--}.
\end{align*}
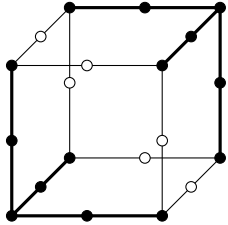
\[
\begin{tikzpicture}[x={(1cm,0cm)},y={(3.85mm, 3.85mm)},z={(0cm,1cm)},baseline=0]
    \draw[very thick] (-1,-1,-1) -- (-1,-1, 1);
    \draw[very thick] (-1,-1,-1) -- (-1, 1,-1);
    \draw[very thick] (-1,-1,-1) -- ( 1,-1,-1);
    \draw (-1,-1, 1) -- (-1, 1, 1);
    \draw (-1,-1, 1) -- ( 1,-1, 1);
    \draw (-1, 1,-1) -- (-1, 1, 1);
    \draw (-1, 1,-1) -- ( 1, 1,-1);
    \draw ( 1,-1,-1) -- ( 1,-1, 1);
    \draw ( 1,-1,-1) -- ( 1, 1,-1);
    \draw[very thick] (-1, 1, 1) -- ( 1, 1, 1);
    \draw[very thick] ( 1,-1, 1) -- ( 1, 1, 1);
    \draw[very thick] ( 1, 1,-1) -- ( 1, 1, 1);
    
    \draw[fill=black] (-1,-1,-1) circle (2pt) node[above right] {};
    \draw[fill=black] (-1,-1, 1) circle (2pt) node[below] {};
    \draw[fill=black] (-1, 1,-1) circle (2pt) node[above] {};
    \draw[fill=black] (-1, 1, 1) circle (2pt) node[left] {};
    \draw[fill=black] ( 1,-1,-1) circle (2pt) node[right] {};
    \draw[fill=black] ( 1,-1, 1) circle (2pt) node[below] {};
    \draw[fill=black] ( 1, 1,-1) circle (2pt) node[above] {};
    \draw[fill=black] ( 1, 1, 1) circle (2pt) node[below right] {};
    
    \draw[fill=black] ( 0,-1,-1) circle (2pt) node[below] {};
    \draw[fill=white] ( 0,-1, 1) circle (2pt) node[below] {};
    \draw[fill=white] ( 0, 1,-1) circle (2pt) node[above] {};
    \draw[fill=black] ( 0, 1, 1) circle (2pt) node[above] {};
    
    \draw[fill=black] (-1, 0,-1) circle (2pt) node[left] {};
    \draw[fill=white] (-1, 0, 1) circle (2pt) node[left] {};
    \draw[fill=white] ( 1, 0,-1) circle (2pt) node[right] {};
    \draw[fill=black] ( 1, 0, 1) circle (2pt) node[right] {};
    
    \draw[fill=black] (-1,-1, 0) circle (2pt) node[below right] {};
    \draw[fill=white] (-1, 1, 0) circle (2pt) node[above left] {};
    \draw[fill=white] ( 1,-1, 0) circle (2pt) node[below right] {};
    \draw[fill=black] ( 1, 1, 0) circle (2pt) node[above left] {};
\end{tikzpicture}
\]

Moreover, we choose \(O_{19} \coloneqq L_{+-0}\) as origin. Then, \(P_{19} \coloneqq L_{0+-}\) and \(Q_{19} \coloneqq L_{-0+}\) become the two \(3\)-torsion sections, because it holds \(\langle P_{19}, P_{19} \rangle = \langle Q_{19}, Q_{19} \rangle = 0\), whereas
\(R_{19} \coloneqq L_{-+0}\) becomes a section of infinite order. From~\cite[eq. (8.12) and Table (8.16)]{Shioda:on.MW.lattices}  it follows
\[
    \langle R_{19}, R_{19} \rangle = 2\chi + 2 L_{-+0}\cdot L_{+-0} -\sum \textit{contr}_v (R_{19})= 2\cdot 2 + 2 \cdot 0 - 2\cdot \frac{4}{3}= \frac{4}{3}.
\]

\begin{theorem} \label{thm:Apery_Fermi}
There is an Enriques involution \(\varepsilon \in \Aut(X)\) which acts on the \(L\)-lines by exchanging all subscripts `\(+\)' with `\(-\)' and on the \(M\)-lines by exchanging \(M_{k+\beta}\) with \(M_{k-\beta}\), for all \(k \in \set{1,2,3}\), \(\beta \in \set{+,-}\). Moreover, \(\varepsilon\) is not of Barth--Peters type.
\end{theorem}
\proof
Let \(S_{19} \coloneqq \boxminus R_{19}\) be the section given by the inverse of \(R_{19}\) in the Mordell--Weil group. 
Then we clearly have \(\langle S_{19}, S_{19}\rangle = \langle R_{19}, R_{19}\rangle\) and \(\langle S_{19}, R_{19}\rangle = -\langle R_{19}, R_{19}\rangle\).
From these equalities, using~\cite[Theorem 8.6]{Shioda:on.MW.lattices} we obtain \(O_{19} \cdot S_{19}=0\) and \(R_{19} \cdot S_{19}=2\).
These intersection numbers explicitly determine \(S_{19}\) in the coordinate system of \autoref{rmk:Hulek.Schutt.II^*}:
\[
S_{19}= (19, 17, -27, -42, -54, -81, -66, -51, -34, -17, -27, -42, -54, -81,
-66, -51, -34, -17, 4)\, .
\]

We are then able to  compute  the translation by \(R_{19}\), denoted by \(t\), and involution \(\imath\) as in Hulek and Schütt's construction \cite[§3]{Hulek.Schuett:Enriques.surf.jacobian.ell}.
Explicit computations show that the invariant lattice of \(\varepsilon \coloneqq t \circ \imath\) is isomorphic to \(\bM\), so that \(\varepsilon\) is the pullback of an Enriques involution. 
We can verify directly that \(\varepsilon\) acts on the \(L\)-lines and \(M\)-lines as described in the statement of the theorem. 
By~\autoref{lem:emb.M-->S_X.for.T_X=U+4m} and the proof of \autoref{thm:Enr(X)_Apery_Fermi},
we know that the coinvariant lattice \((S_X)_\varepsilon\) is isomorphic to either \(\bA_2(2) \oplus \bE_7(2)\) or \( \bE_8(2) \oplus [-12] \).
An explicit computation shows that \((S_X)_\varepsilon\) contains \(132\) vectors of square \(-4\), so it is necessarily isomorphic to \(\bA_2(2) \oplus \bE_7(2)\), i.e. \(\varepsilon\) is not of Barth--Peters type.
We refer to the ancillary file \texttt{calc.Apery\char`_Fermi.sage} for the actual computations in Sage.
\endproof

\begin{remark} \label{rmk:Apery-Fermi.ell.fibr}
Thanks to the description of the Enriques involution in \autoref{thm:Apery_Fermi}, it is immediate to see that the images of the \(L\)-lines in \(Y\) form a tetrahedron, while the images of the \(M\)-lines form a complete graph with \(6\) vertices in which three pairs of curves intersect doubly. 
The tetrahedron and the complete graph are connected in the following way, where double intersections are marked with a double edge.
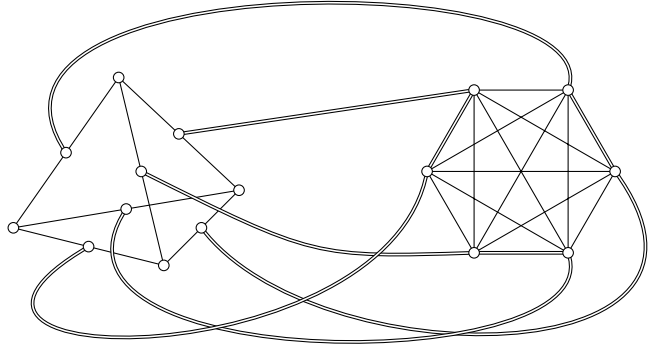
\[
\begin{tikzpicture}
\coordinate (A) at (0,0);
\coordinate (B) at (2,-0.5);
\coordinate (C) at (3,0.5);
\coordinate (D) at (1.4,2);

\coordinate (Z) at (8,0.75);

\draw (A) -- (D) coordinate [midway] (E);
\draw (B) -- (C) coordinate [midway] (F);
\draw (A) -- (B) coordinate [midway] (G);
\draw (C) -- (D) coordinate [midway] (H);
\draw (A) -- (C) coordinate [midway] (I);
\draw (B) -- (D) coordinate [midway] (L);

\path (Z) arc (0:0:1.25) coordinate (M) 
      (Z) arc (0:60:1.25) coordinate (N) 
      (Z) arc (0:120:1.25) coordinate (O)
      (Z) arc (0:180:1.25) coordinate (P)
      (Z) arc (0:240:1.25) coordinate (Q)
      (Z) arc (0:300:1.25) coordinate (R);

\draw[double] (M) -- (N);
\draw (M) -- (O) (M) -- (P) (M) -- (Q) (M) -- (R) (N) -- (O) (N) -- (P) (N) -- (Q) (N) -- (R);
\draw[double] (O) -- (P);
\draw (O) -- (Q) (O) -- (R) (P) -- (Q) (P) -- (R);
\draw[double] (Q) -- (R);

\draw[double] (H) -- (O);
\draw[double] (L) .. controls (4,-0.4) .. (Q);
\draw[double] (F) .. controls (4,-2) and (10,-2) .. (M);
\draw[double] (E) .. controls (-1, 3.5) and (8, 3.5) .. (N);
\draw[double] (I) .. controls (0,-2) and (8,-2) .. (R);
\draw[double] (G) .. controls (-2,-2) and (5,-2) .. (P);

\foreach \node in {A,B,C,D,E,F,G,H,I,L,M,N,O,P,Q,R}
    \draw[fill=white] (\node) circle (2pt);
\end{tikzpicture}
\]

The following pencils (we omit here the images of the \(M\)-lines) are special pencils on \(Y\) whose pullbacks correspond to the elliptic fibrations on \(X\) with frames \(W_{19}\), \(W_{20}\), \(W_{24}\) and \(W_{26}\), respectively. 

\[
\begin{tikzpicture}
\coordinate (A) at (0,0);
\coordinate (B) at (2,-0.5);
\coordinate (C) at (3,0.5);
\coordinate (D) at (1.4,2);

\draw[very thick] (A) -- (D) coordinate [midway] (E);
\draw (B) -- (C) coordinate [midway] (F);
\draw (A) -- (B) coordinate [midway] (G);
\draw[very thick] (C) -- (D) coordinate [midway] (H);
\draw (A) -- (C) coordinate [midway] (I);
\draw[very thick] (B) -- (D) coordinate [midway] (L);

\foreach \node in {A,B,C,D,E,L,H}
    \draw[fill=black] (\node) circle (2pt);

\foreach \node in {I,G,F}
    \draw[fill=white] (\node) circle (2pt);
\end{tikzpicture}
\quad
\begin{tikzpicture}
\coordinate (A) at (0,0);
\coordinate (B) at (2,-0.5);
\coordinate (C) at (3,0.5);
\coordinate (D) at (1.4,2);

\draw (A) -- (D) coordinate [midway] (E);
\path (B) -- (C) coordinate [midway] (F);
\path (A) -- (B) coordinate [midway] (G);
\draw (C) -- (D) coordinate [midway] (H);
\path (A) -- (C) coordinate [midway] (I);
\draw (B) -- (D) coordinate [midway] (L);

\draw[very thick,densely dashed] (A) -- (B);
\draw[very thick,densely dashed] (B) -- (C);
\draw[very thick,densely dashed] (C) -- (A);

%\draw[thick] [yshift=-1pt] (0,0) -- (2,-0.5) [yshift=2pt] (0,0) -- (2,-0.5);
%\draw[thick] [yshift=-1pt] (0,0) -- (3,0.5) [yshift=2pt] (0,0) -- (3,0.5);
%\draw[thick] [yshift=-1.5pt] (2,-0.5) -- (3,0.5) [yshift=3pt] (2,-0.5) -- (3,0.5);

\foreach \node in {A,G,B,F,C,I,D}
    \draw[fill=black] (\node) circle (2pt);

\foreach \node in {E,L,H}
    \draw[fill=white] (\node) circle (2pt);
\end{tikzpicture}
\quad
\begin{tikzpicture}
\coordinate (A) at (0,0);
\coordinate (B) at (2,-0.5);
\coordinate (C) at (3,0.5);
\coordinate (D) at (1.4,2);

\draw (A) -- (D) coordinate [midway] (E);
\draw (B) -- (C) coordinate [midway] (F);
\draw (A) -- (B) coordinate [midway] (G);
\draw (C) -- (D) coordinate [midway] (H);
\draw[very thick] (A) -- (C) coordinate [midway] (I);
\draw (B) -- (D) coordinate [midway] (L);

\draw[very thick] (E) -- (A) -- (G);
\draw[very thick] (H) -- (C) -- (F);

\foreach \node in {A,C,E,F,G,H,I,L}
    \draw[fill=black] (\node) circle (2pt);

\foreach \node in {B,D}
    \draw[fill=white] (\node) circle (2pt);
\end{tikzpicture}
\quad
\begin{tikzpicture}
\coordinate (A) at (0,0);
\coordinate (B) at (2,-0.5);
\coordinate (C) at (3,0.5);
\coordinate (D) at (1.4,2);

\draw[very thick] (A) -- (D) coordinate [midway] (E);
\draw[very thick] (B) -- (C) coordinate [midway] (F);
\draw[very thick] (A) -- (B) coordinate [midway] (G);
\draw[very thick] (C) -- (D) coordinate [midway] (H);
\draw (A) -- (C) coordinate [midway] (I);
\draw (B) -- (D) coordinate [midway] (L);

\foreach \node in {A,B,C,D,E,F,G,H}
    \draw[fill=black] (\node) circle (2pt);

\foreach \node in {I,L}
    \draw[fill=white] (\node) circle (2pt);
\end{tikzpicture}
\]
\end{remark}

\bibliographystyle{amsplain}
\bibliography{references}

\end{document}